\input amstex
\loadbold
\documentstyle{amsppt}

\NoBlackBoxes
\magnification=1100
\baselineskip19pt

\define\Ob{\operatorname{Ob}}
\define\Hom{\operatorname{Hom}}
\define\ct{T^*X}
\define\cty{T^*Y}
\define\ode{\operatorname {E}}

\define\od{\operatorname  {D}}
\define\op{\operatorname {P}}
\define\om{\operatorname {M}}
\define\oc{\operatorname {C}}
\define\oss{\operatorname {Char}}     

\define\cross{\times}
\define\intersect{\cap}
\define\union{\cup}
\define\bigunion{\bigcup}  
\define\muhom{\mu \boldsymbol h \boldsymbol o \boldsymbol m}
\define\var{\operatorname{var}}

\define\RE{\operatorname{Re}}        

\topmatter
\title Microlocal perverse sheaves 
\endtitle
\author S. Gelfand, R. MacPherson and K. Vilonen
\endauthor
\address American Mathematical Society, 201 Charles street, Providence, RI  02904
\endaddress
\email SXG\@math.ams.org
\endemail
\address Institute for Advanced Study, Princeton NJ 08540
\endaddress
\email rdm\@math.ias.edu
\endemail
\address{Department of mathematics, Northwestern University, Evanston IL 60208}
\endaddress
\email Vilonen\@math.northwestern.edu
\endemail
\thanks We thank the I.H.E.S. for its hospitality in enabling us to begin our  work on this
 subject
\endgraf  K.Vilonen was partially supported by NSA, NSF, and DARPA.  
\endthanks

\abstract{\noindent Microlocal perverse sheaves form a stack on the cotangent bundle $\ct$ of a complex manifold
that is the analogue of the stack of perverse sheaves on the manifold $X$ itself.   We give an embedding of the stack
of microlocal perverse sheaves into a simpler stack on $\ct$  which we call {\it melded Morse systems}.  The
category of melded Morse systems is elementary, in the sense that an object is a collection of vector spaces and a
collection of linear maps between them satisfying certain composition relations.}\endabstract
\endtopmatter

\document

\bigpagebreak

\heading{0. Introduction.}
\endheading

In this paper we study microlocal perverse sheaves. The category of perverse sheaves on a complex manifold $X$ is
equivalent to the category of holonomic regular $\Cal D_X$-modules on $X$; here $\Cal D_X$ stands for the sheaf of linear
differenatial operators on $X$. In the same manner the category of microlocal perverse sheaves is
equivalent to the category of holonomic regular $\Cal E_X$ modules on $\ct$; here $\Cal E_X$ denotes the sheaf of
microdifferential operators on $\ct$. 

To any open set $U$ in a complex algebraic variety $X$, we have the
category $\op(U)$ of perverse sheaves on $U$.  Perverse sheaves have many beautiful formal properties.  Two of
the most important are: 

\roster
\item $\op(U)$ is an abelian category.

\item The assignment $U \mapsto \op (U)$ is a stack, or ``sheaf of categories".  See the appendix for a
summary of what we need from the theory of stacks.  

\endroster

The category $\op (U)$ is usually constructed as a subcategory of $D^b(U)$, the constructible derived category of
the category of sheaves on $U$ (see \cite{BBD}).  This construction has the  aesthetic deficiency that both of the
two properties listed above fail for the ambient category: $D^b(U)$ is not abelian, and the assignment $U \mapsto
D^b(U)$ is not a stack.   In this paper, we will remedy this deficiency by microlocalization, as will be explained
below.  

There is another stack of abelian categories $\ode$ on
$\ct$, called the stack of microlocal perverse sheaves.  This is a refinement of the stack $\op$ on $X$ in the
following sense:  There is an equivalence of categories called the {\it microlocalization functor} $\mu$ from the
category $\op (U)$ of perverse sheaves on U to the category
$\ode (\pi^{-1} (U))$ of microlocal perverse sheaves on the preimage of $U$ in $\ct$.   For an open set $V \subset
\ct$, the category
$\ode (V)$ is constructed as a full subcategory Kashiwara and Schapira's microlocal derived category $D^b(X,V)$
[KS] (see \S 2 below).   But like $D^b(U)$, the ambient category $D^b(X,U)$ fails to be abelian and doesn't form a
stack.

 In this paper, we will construct an embedding $\nu: \ode \to \om$ of
the stack $\ode$ of microlocal perverse sheaves into a new stack $\om$ of melded Morse systems.  The construction
of $\om$ involves no derived categories; it is purely geometric (although somewhat technical).  The category $\om
(V)$ is visibly abelian, and the assignment $V
\mapsto \om(V)$ is clearly a stack.  

In order give the idea behind $\om$, let's consider the composed functor $\nu \circ \mu: \op 
\to \ode \to \om$.  Given a perverse sheaf $\Cal P$, we will explain what kind of data this composed functor will
assign to it.  First of all, associated to $\Cal P$ there is a Lagrangian variety $\Lambda\subset\ct$ called the { \it characteristic
variety} of $\Cal P$, which is the support of $\mu (\Cal P)$.  

For simplicity in the rest of this introduction we will assume that $X$ is endowed with a fixed Whitney
stratification $\Cal S$.  This is a disjoint decomposition  $X = \bigcup_\alpha \Cal S_\alpha$ of $X$ into strata
$\Cal S_\alpha$ that are smooth connected algebraic subvarieties.  We will consider only perverse sheaves that are
constructible with respect to the stratification $\Cal S$.  Then $\Lambda = \bigcup_\alpha T^*_{\Cal S_\alpha}X$ is the
union of the conormal bundles $T^*_{\Cal S_\alpha}X$ to the strata $\Cal S_\alpha$.  This is a closed singular Lagrangian
subvariety of $\ct$.  The fact that a perverse sheaf is constructible with respect to the stratification $\Cal S$ implies that its
microlocalization has support on $\Lambda$ (and conversely). 

The variety $\Lambda \subset \ct$ itself has a stratification.  We define $\Lambda^0$ to be the nonsingular part of
$\Lambda$.  It is the union of connected $n$ dimensional manifolds each of the form $\Lambda^0 \cap T^*_{\Cal
S_\alpha}X$, which is the part of
$T^*_{\Cal S_\alpha}X$ that is not in the closure of the conormal bundle $T^*_{\Cal S_\beta}X$ of any stratum
$\Cal S_\beta$, for
$\beta \neq \alpha$.)  We denote by $\Lambda^1$ the union of the codimension one strata (see \S 4, Step 1), so that
$\Lambda ^2 =
\Lambda - \Lambda^0 - \Lambda ^1$ has codimension two in $\Lambda$. The space $\Lambda^2$ is itself the
union of strata of various dimensions, but we will not need these in this paper since we have not been able to
understand microlocal perverse sheaves on
$\Lambda^2$.  

In this paper we will describe the structure of the stack of microlocal perverse sheaves on the open sets
$(\ct)^0 = \ct - \Lambda^1-\Lambda^2$ and $(\ct)^1 = \ct - \Lambda^2$. We do so by applying  complex  ``Morse"
theory to the perverse sheaf $\Cal P$.  The {\it critical points} of a proper complex analytic function $f: X \to \Bbb C$ are the
points $x \in X$ such that $df(x) \in \ct$  is in the characteristic variety of $\Cal P$, and the {\it critical values} are the values
$v \in \Bbb C$ such that $v = f(x)$ for some critical point $x$.  The sheaf $f_* (\Cal P)$ is singular at the critical values.
To define $\om$ on $(\ct)^0$ we consider all Morse functions on $(\ct)^0$, i.e., functions $f$ at points $x$ such that $df$ is
transverse to $\Lambda$ and meets $\Lambda$ at a single point in $df(x)\in\Lambda^0$. Then $\om$ on $(\ct)^0$ consists of
all the assginments $f\mapsto f_* (\Cal P)$; here we must view $f_* (\Cal P)$ in localized category as is explained in \S 3.
In this case there is only one critical value. To extend this definition to $(\ct)^0$ we must explain what happens at points of
$\Lambda^1$. Here we simply consider all functions $f$ such that $df$ intersects $\Lambda$ at exactly on point which lies in
$\Lambda^1$. We furthermore assume that $df$ is in as generic a position as possible with respect to $\Lambda$ at the point
of the intersection. Then we slightly perturb $f$ and again make the assignment $f\mapsto f_* (\Cal P)$, as before. Now,
contrary to the case of $\Lambda^0$, we will in general have many critical values. 

As was mentioned earlier, we have not been able to analyse the conditions that need to be imposed on objects in $\om$
in order for them to lie in the image of $\ode$. It is generally belived (a conjecture of Kashiwara) that no new conditions
need to be imposed beyond the codimension 2 locus $\Lambda^2$. 

This paper is an updated version of a manuscipt that goes back a few years. 

We wish to thank A. Beilinson, T. Braden, M. Kashiwara, and D. Nadler for conversations. The conversations with
Kashiwara were especially helpful for some of the  technical details of the paper.

\heading 1. Symplectic manifolds, contact 
 transformations and the Maslov bundle.
\endheading

This section contains a review of symplectic geometry most of which can be found, for example, in the appendix of
\cite{KS}. Let $M$ be a complex manifold. Recall that if
$\omega$ is a non-degenerate closed 2-form  on $M$, then $(M,\omega)$ is called a symplectic manifold. Via the
two form $\omega$, we can identify (locally defined) vector fields and (locally defined) 1-forms as follows. The
vector field $v$ corresponds to the  1- form $\alpha$ if  $$ 
\alpha (w) = \omega (w,v)
$$ for all vector fields $w$.

Let us now assume that $M$ has a $\Bbb C^*$-action. This $\Bbb C^*$-action gives rise, by differentiating, to a
vector field $e$ on
$M$. Let $\alpha$ be the  1-form corresponding to $e$. If $\omega = d\alpha$ we call $M$ a  homogenous
symplectic manifold.  

Let $X$ be a complex manifold. Then the conormal bundle $\ct$ is naturally a homogenous  symplectic manifold as
follows. Let $\pi :
\ct\to X$ be the projection. We define a 1-form $\alpha$ on $\ct$ by the formula $\alpha (v) = \xi (\pi (v))$, where
$\xi
\in \ct$ and $v \in T_\xi(\ct)$. Then $\omega = d\alpha$ gives a symplectic structure on $\ct$.  Since $\ct$ is a vector
bundle over
$X$, it has a natural $\Bbb C^*$-action. This makes $\ct$ into a homogenous symplectic manifold.

The easiest way to see this is to use local coordinates. Let $(x_1,\dots,x_n)$ be local coordinates on $X$. We get
local coordinates ($x_1,\dots,x_n,\xi_1,\dots, \xi_n)$ on
$\ct$ by writing $\xi\in T_x^*X$, $x=(x_1,\dots,x_n)$ as $\xi = \sum\xi_i dx_i$. Then, expressed with respect to
these local coordinates on $\ct$,  $\alpha = \sum \xi_idx_i$,
$e=\sum  x_i\dfrac {\partial}{\partial x_i}$ and $\omega = \sum d\xi_i\wedge dx_i$.

An analytic subvariety $\Lambda$ of a symplectic manifold $M$ is called {\sl Lagrangian} if at any smooth point
$\xi\in\Lambda$ we have $T_\xi\Lambda = (T_\xi\Lambda)^\perp$, where $(T_\xi\Lambda)^\perp$ is the
orthogonal  complement of $T_\xi\Lambda$ in $T_\xi M$ with respect to $\omega$. If $M$ is a homogenous
symplectic manifold then a subset is called {\sl conic} if it is invariant under the
$\Bbb C^*$-action. Standard examples of conic Lagrangian subvarieties in $\ct$ are conormal bundles 
$$ T_S^*(X)=\{(x,\xi)\in\ct\mid x\in S, \xi|T_x(S)=0\}
$$ to smooth subvarieties $S$ in $X$.

Let $\Lambda \subset \ct$ be a conic Lagrangian subvariety. We say that $\Lambda$ is {\sl in generic position} at
$\xi\in\Lambda$ if
$\Lambda\cap\pi^{-1}(\pi(\xi)) $ is a one dimensional complex subspace of the fiber
$\pi^{-1}(\pi(\xi)) $. Let us note that $T_S^*(X)$ is in generic position if and only if $\operatorname{codim} S =1$.

A standard technique in symplectic geometry is to reduce the study of arbitrary conic Lagrangian subvarieties
$\Lambda \subset \ct$ to the study of those in generic position.  This is accomplished in two steps.  The first step is
to get rid of the zero section as a possible component of $\Lambda$.  This is accomplished by embedding $X$ as a
codimension one subvariety of another variety, say
$X=X\times \{0\}\subset X\times
\Bbb C$, and replacing $\Lambda \subset \ct$ by the Lagrangian subvariety $\Lambda
\times T_{\{0\}}^* \Bbb C \subset T^*(X\times \Bbb C)=T^*X\times T^*\Bbb C$.  The second step is to apply a
contact transformation, as described in the proposition below.

Let $U\subset\ct$ and $V\subset\cty$ be conic, open subsets. Let $\omega_X$ and
$\omega_Y$ be the symplectic forms and $\alpha_X$ and $\alpha_Y$ be the canonical 1-forms on $\ct$ and $\cty$.
A biholomorphic map
$\phi : U\to V$ is called a {\sl contact transformation} if it is $\Bbb C^*$-equivariant and $\phi^*\omega_Y =
\omega_X$.  Equivalently, $\phi$ is a contact  transformation if $\phi^*\alpha_Y =  \alpha_X$.

Let $\phi : U\to V$ be a biholomorphic map and let $\Gamma\subset T^*(X\times Y)$ be the graph of $\phi$. Then
$\phi$ is a contact transformation if and only  if $\Gamma$ is Lagrangian with respect to the symplectic structure
$(\omega_X,-\omega_Y)$ on
$T^*(X\times Y)$. We call $\phi$ a  generic contact  transformation if $\Gamma$ is the intersection of $U\times V$
with $ T^*_S(X\times Y)$, for some smooth hypersurface
$S\subset X\times Y$. 

\proclaim {Proposition 1.1} Let $\Lambda\subset\ct$ be a Lagrangian subvariety in which the zero section does not
appear as a component,  and let $\xi\in\Lambda$. Then there exists a neighborhood $U$ of $\xi$ and a  generic
contact transformation $\phi : U\to V\subset\cty$ such that $\phi (U\cap\Lambda)$ is in generic position.
\endproclaim

\demo{Proof} See \cite{KK, Corollary 1.6.4}.
\enddemo

We will now define the complex analogue of the Maslov bundle. Let
$(V, \omega)$ be a complex symplectic vector space of dimension $2n$ and let  $Gr(V)$ be  the Lagrangian
Grassmannian, i.e., the space of all Lagrangian vector subspaces of
$V$.

Let us fix three integers $n_{12}$, $n_{23}$, and $n_{13}$ such that $0\leq n_{ij} \leq n$. Consider the following
variety
$$ X=\{(W_1,W_2,W_3)\in Gr(V)\times Gr(V)\times Gr(V) \ | \  
\text{dim}(W_i\cap W_j)= n_{ij}\}. $$

On $X$ we have the tautological bundle $\Cal W \subset X\times V\times V \times V$ whose fibre over
$(W_1,W_2,W_3)$ is $W_1\times W_2\times W_3\subset V\times V\times V$. On $\Cal W$ we have a quadratic
form $Q$ defined as follows. For
$(x_1,x_2,x_3)\in W_1\times W_2\times W_3$ define
$$ Q(x_1,x_2,x_3)= \omega (x_1,x_2) + \omega (x_2,x_3) + \omega (x_3,x_1). $$

On each fibre the quadratic form $Q$ has rank $m=3n-n_{12}-n_{23}- n_{13}$.

We now consider $\operatorname{Re} Q$ and let $\Cal W' \subset \Cal W$  be a  continuous real subbundle such
that on each fibre
$W\subset  W_1\times  W_2\times W_3$ of $\Cal W'$ the form $\operatorname{Re} Q$ is negative definite and
$W$ is maximal with respect to this property.

\proclaim {Lemma 1.2} The subbundles $\Cal W'$
 exist and they are all isomorphic.
\endproclaim

\demo {Proof} Let $E\to X$ be the bundle of real
$(3n-n_{12}-n_{23}-n_{13})$-dimensional subspaces of $\Cal W$. Denote by $E'$ subbundle of $E\to X$
consisting of subspaces $W$  such that $\operatorname{Re} Q$ is negative definite on $W$. We compute the fibre
of $E'$ at a point 
$(W_1,W_2,W_3)\in X$.

On the vector space $W_1 \oplus W_2 \oplus W_3$ we choose coordinates  such that the quadratic form
$\operatorname{Re} Q$ is given by the matrix $$
\left(
\matrix -\operatorname I & 0 & 0 \\ 0 & \operatorname I & 0  \\ 0 & 0 & 0 
\endmatrix
\right).
$$

Let $W$ be the subspace corresponding to the block $-\operatorname{I}$, $W^\perp$ the subspace corresponding
to the block
$\operatorname {I}$ and let $N$ be the nullspace of $\operatorname{Re} Q$. We have $\dim_{\Bbb R} W =
\dim_{\Bbb R} W^\perp = 3n-n_{12}-n_{13}-n_{23}$, $\dim_{\Bbb R} N= 2(n_{12}+n_{13} + n_{23})$. If
$W'$ is any other $3n-n_{12}-n_{13} -n_{23}$-dimensional subspace on which $\RE Q$ is negative definite, then
$W'\subset W\oplus W^\perp\oplus N$ and
$W'$ projects onto $W$.

Therefore $W'$ is given by a graph of a linear function $f:W\to W^\perp \oplus N$ and we write $f(w) = f_1(w) +
f_2(w)$ with
$f_1(w)\in W^\perp$ and $f_2(w)\in N$. For $W'$ to be negative definite we must have
$$ 0 > \RE Q(w+f(w)) = -|w|^2 + |f_1(w)|^2.
$$

This means precisely that we must have $|f_1| < 1$. Therefore the bundle $E'$ has contractible fibres. Because the
bundles $\Cal W'$ are given by sections of $E'$ they exist and are all isomorphic.
\enddemo

Given such a bundle $\Cal W'$ we form $\Lambda^m\Cal W'$. It is a real line bundle on
$X$, and therefore gives rise to an element in $H^1(X,\Bbb Z/2\Bbb Z)$.

We therefore get a local system $\mu$ with structure group $\{\pm 1\}$. We call this local system $\mu$ the {\sl
universal Maslov bundle.}

\heading 2. Microlocal perverse sheaves
\endheading

In this section we define microlocal perverse sheaves as a stack $\ode$ on the cotangent bundle $\ct$.  (In the
appendix, the definition and properties of a stack are recalled.)    For
$U\subset \ct$, we define the microlocal derived category $D^b(X,U)$ by algebraically localizing the derived
category of sheaves on
$X$ away from $U$. The association $U \mapsto D^b(X,U)$ is a prestack.   Microlocal perverse sheaves are cut out
of this prestack by appropriate constructibility and perversity restrictions.

We will work with conic open sets in $\ct$.  The set of all conic open sets in $\ct$ is a (non-Hausdorff) topology on
$\ct$, the {\it conic} topology.  The stack $\ode$ will be a stack with respect to the conic topology.  The identity
map from $\ct$ with the usual topology to $\ct$ with the conic topology is continuous.  Therefore, by the pullback
construction on stacks, we also get a stack with the usual topology, but it contains no more information.

Let $X$ be a complex manifold and let $\op_X$ denote the stack of perverse  sheaves of $\Bbb C$-vector spaces on
$X$. We will use the conventions of \cite{BBD} in indexing the complexes. Given a perverse sheaf $A \in
\op_X(W)$ on an open subset $W \subset X$ we have the characteristic variety (or singular support) $\oss(A)
\subset T^*U$, which is a conic Lagrangian variety, whose definition is recalled below. 

Let $\od^b(X)$ denote the bounded derived category of $\Bbb C$-sheaves on $X$, and let $\od^b_{\Bbb R-c}(X)$
denote the subcategory of $\od^b(X)$ consisting of complexes whose cohomology is  constructible with respect to a
real subanalytic stratification.

For any $A\in\od^b(X)$ we define its singular support, or characteristic variety, $\oss (A)
\subset \ct$ following \cite{KS} as follows. We say that $\xi\notin\oss (A)$ if there exists an open set
$U\ni\xi$ such that for any $x\in X$ and for  any smooth function $f$ defined in a neighborhood of $x$ with
$f(x)=0$ and $df_x\in U$ we have  $$   
\operatorname{R\Gamma}_Z (A)_x=0,
$$ where $Z=\{ y\in X | f(y) \geq 0\}$.

This definition immediately implies that $\oss (A)$ is a closed $\Bbb R^+$-conic subset of $\ct$. Let $U \subset \ct$
be a $\Bbb C^*$-conic (this is a  hypothesis needed in the future and not necessary at the moment) open  subset.  As
in \cite{KS}, we consider the following subset $\Cal N$ of objects in $\od^b(X)$, $\Cal N = \{ A\in \od^b(X) |  \oss
(A) \cap U = \emptyset \}$. This subset satisfies the following conditions

\roster
\item $0\in\Cal N$.
\item $A\in \Cal N \Rightarrow A[1] \in \Cal N$.
\item If $A\to B \to C$ is a distinguished
 triangle then  $B \in \Cal N$ whenever $A\in \Cal N$ and $C\in \Cal N$.
\endroster

We can now localize $\od^b(X)$ with respect to $\Cal N$ (see, e.g., \cite{KS, Chap. I,
\S 6}) to get a category $\od^b(X,U)$. In this category all the objects in $\Cal N$ become isomorphic to the zero
object. We also have a natural functor  $\od^b(X)  \to
\od^b(X,U)$. Note that $U \mapsto \od^b(X,U)$ is a  prestack on $\ct$. For $A\in
\od^b(X,U)$ the set $\oss (A) \cap U$ is well defined. In a similar way, we can localize 
$\od^b_{{\Bbb R}-c}(X)$ to  form the localized categories $\od^b_{\Bbb R-c}(X,U)$.

  Let $U\subset \ct$, $V\subset \cty$, and let $\phi : U \to V$ be a contact  transformation which we assume to be
generic. Then 
$\text{graph}(\phi) =T^*_S(X\times Y)$, where $S\subset X\times Y$ is a smooth complex hypersurface. We get
two functors $\theta\:
\od^b(X) \to \od^b(Y)$ and $\tau\: \od^b(Y) \to
\od^b(X)$ given by 
$$
\align
\theta  &= Rq_*p^*[n-1], \\
\tau  &= Rp_*q^* [n-1],
\endalign
$$ where $p: S\to X$ and $q : S \to Y$ are the projections and $n=\dim X = \dim Y$. If
$S$ is subanalytic and $\overline S$ is compact, then the functors $\theta$ and $\tau$ preserve constructibility. Note
that by making the open sets $U$ and $V$ smaller it is possible to make $S$ satisfy these hypotheses.

\proclaim {Proposition 2.1} The functors $\theta$ and $\tau$ descend to  functors
$\overline\theta : \od^b(X,U) \to \od^b(Y,V)$ and $\overline\tau : \od^b(Y,V) \to
\od^b(X,U)$, which are equivalences of categories inverse to each other. For
$A\in\od^b(X,U)$ we have $\phi(\oss (A) \cap U) = \oss (\overline\theta  (A))$. If $S$ is subanalytic and $\overline
S$ is compact then the statement holds for the constructible  categories $\od^b_{\bold R-c}(X,U)$ and 
$\od^b_{\bold R-c}(Y,V)$. 
\endproclaim

\demo {Proof} See \cite{KS, Chap. VII}.
\enddemo

Recall that the characteristic variety of a sheaf in $\od^b_{\Bbb R-c}(X)$ is a subanalytic, $\Bbb R^+$ conic, real
Lagrangian subvariety of $T^*(X)$.   Let $U\subset\ct$ be a $\Bbb C^*$-conic open subset of $\ct$. We say that
$A\in\od^b_{{\bold R}-c}(X)$ is microlocally complex constructible in $U$ if $\Lambda = \oss (A) \cap U$ is a
complex conic Lagrangian variety.

\proclaim {Proposition 2.2} If $A$ is microlocally complex constructible in $\ct$, then it is complex constructible in
the usual sense. \endproclaim

\demo {Proof} Follows from Theorem 8.5.5 in \cite{KS}.
\enddemo

\subheading {Preliminaries on Morse functions} Here we define a class of Lagrangian manifolds ${\Cal T} \subset
T^*X$, called test manifolds. 

Let $U\subset\ct$ be a $\Bbb C^*$-conic open subset and let $\Lambda\subset U$ be a closed conic Lagrangian
subvariety. Let
$X=\bigunion X_\alpha$ be a subanalytic stratification of $X$ such that $\Lambda \subseteq \bigunion_\alpha
T^*_{X_\alpha}X$.  We consider triples $(B,f,D)$, where $B\subset X$ is a closed differentiably embedded ball,
$f:U\to\Bbb C$ is a complex analytic function defined on some open neighborhood of $B$ in $X$, and $D\subset
f(B) \subset \Bbb C$ is a closed differentiably embedded disk. The triple  $(B,f,D)$ will be called a {\sl test triple}
if there is a neighborhood
$V\subset\Bbb C$ of $D$ such that for all strata $X_\alpha$, the differential $df$ never vanishes on $\partial
B\intersect X_\alpha
\intersect f^{-1}(V)$, where $\partial B$ is the  boundary of $B$. It follows from Thom's first isotopy theorem that
$\partial B\intersect f^{-1}D$ is homeomorphic to the product of $D$ with a stratified space and $f$ is equivalent to
the projection to $D$.  

To each such triple $(B,f,D)$, we can associate a complex  analytic Lagrangian  submanifold ${\Cal T}\subset
T^*X$ by ${\Cal T}=\{df(x)|x\in B^0\intersect f^{-1}D^0\}$, where $B^0$ and $D^0$ are the interiors of $B$ and
$D$.  We say that
${\Cal T}$ is a {\sl test manifold} if ${\Cal T}$ arises in this way  from a test triple
$(B,f,D)$. Given ${\Cal T}$, there are several triples $(B,f,D)$ to which it is associated. (For example, $f$ is
determined by ${\Cal T}$ only up to an additive constant, and $B$ is determined only in its intersection with
$f^{-1}(D)$.) We will use the ideas of test manifold and test triple interchangeably, according to convenience.  We
use the notation from the introduction ${\Cal B}=\pi({\Cal T})=B^0\intersect f^{-1}(D^0)$, $h=f\circ
\pi\: {\Cal T}\to\Bbb C$, and we denote the restriction $f|{\Cal B}$ also by $f$.

The point of the definition is this.  Suppose $P$ is a constructible sheaf on $X$  whose singular support is contained
in $\Lambda$, and let $i:{\Cal B}\to X$ be the  inclusion. If $(B,f,D)$ is a test triple, i.e., if ${\Cal T}$ is a test
manifold, then the sheaf 
$Rf_*i^*P$ is constructible with respect to the stratification of $\Bbb C$ by the  finite set $h({\Cal T}\intersect
\Lambda)$ (i.e., there are no singularities  coming from the boundary).

\subheading{Definition of microlocal perverse sheaves} Note that if $A\in\od^b_{\Bbb R-c}(X)$ is microlocally
complex  constructible in $U$ then we can apply to it the complex Morse theory of \cite{GM1} and \cite{GM2}. In
particular, we can define the notion of a microlocal  perverse sheaf. Let $A\in\od^b_{{\bold R}-c}(X)$. We say that
$A$ is {\it microlocally perverse } in a
$\Bbb C^*$-conic open subset $U$ of $\ct$ if it is microlocally complex  constructible in
$U$ and 
$$ \operatorname{H}^i(B,f^{-1}(\epsilon);A) = 0 \ \ \ \text{for} \ \ i\neq 0 \ \ \text{and}\ \  \epsilon \in D-\{0\} $$ 
for any test triple $(B,f,D)$ such that $df(B) \cap \oss(A) \cap U  = \{\xi\}$ is a single point whose projection to $X$
is $x$ and $f(x)=0$

\proclaim {Definition 2.3}  For any $\Bbb C^*$-conic, open subset $U$ of $\ct$ let 
$\op(X,U)$ be the full subcategory of $\od^b_{\Bbb R-c}(X,U)$ consisting of sheaves which are  microlocally
perverse in $U$.  The association of the category $\op(X,U)$ to each $\Bbb C^*$-conic, open subset $U$ of $\ct$
forms an abelian prestack on $\ct$.  The stack of {\it microlocal perverse sheaves} is the sackification of this prestack,
as defined in \S 5 of the Appendix.  We denote it by 
$\ode$.
 \endproclaim

\proclaim{Proposition 2.4}  Let $\phi:U \to V$ be a contact  transformation whose defining hypersurace $S$ is
subanalytic with compact closure.  Then the categories  $\op(X,U)$ and
$\op(Y,V)$ are canonically equivalent under $\bar\phi$. 
\endproclaim This follows from Proposition 2.3 together with the remark that there are enough test triples which
remain test triples after the contact transformation to cut out the category $\op(Y,V)$.

Suppose $\Lambda$ is a conic complex Lagrangian subvariety of $\ct$.  We define $\ode_\Lambda$ to be the full
substack of $\ode$ whose objects are supported on $\Lambda$.  We define $\op_\Lambda$ to be the full substack of
the stack of perverse sheaves $\op$ whose objects have characteristic variety lying in $\Lambda$. We have

\proclaim {Lemma 2.5} Suppose that $\xi$ is a point in $T^*_x (X)$ and the conic Lagrangian subvariety
$\Lambda \subset \ct$ is in general position near $\xi$.  Then
$$ (E_\Lambda)_\xi = (\op_{\Lambda \union X}/\Cal N)_x$$ where $\Lambda \union X$ is the union of
$\Lambda$ with the zero section of
$\ct$, $/ \Cal N$ means localization with respect to the local systems on $X$, and the subscripts $\xi$ resp. $x$
mean taking the stalk of the stacks at $\xi$ resp. $x$.
\endproclaim  

This is proved in \cite{W1} and stated as Proposition 2.1.6 in \cite{W2}.  

Sato, Kashiwara, and Kawai have defined a sheaf $\Cal E_X$ of {\it microdifferential
operators} on the cotangent bundle of $X$, which is an analogue of the sheaf $\Cal D_X$ on $X$.  The following
proposition is a microlocal analogue of the Riemann-Hilbert correspondence.

\proclaim{Proposition 2.6} (Andronikov \cite{A}, Waschkies \cite{W2}) The stack $\ode$ is
equivalent to the stack of regular holonomic $\Cal E_X$-modules. 
\endproclaim

\noindent {\bf Remarks on Lemma 2.5 and Proposition 2.6.}  

We had intended to provide proofs for Lemma 2.5 and Proposition 2.6. However, since the first version of this paper was
written, the papers \cite{W1,W2} appeared and we decided to simply refer to them. As the stack regular holonomic $\Cal
E_X$-modules satisfies an anlogue of Lemma 2.5, one needs to define, in a natural fashion, a functor from regular
holonomic $\Cal E_X$-modules to $\ode$. Then Proposition 2.6 follows from the usual Riemann-Hilbert correspondence.
One way to construct such a functor is explained in \cite{W2}. As to Lemma 2.5, it can be proved by a cut-off technique.

\heading 3. Families of perverse sheaves on a disk
\endheading

In this section we introduce, following \cite{GMV}, the technical notion of localized families of perverse sheaves
on a disk. This notion will be used to define melded Morse systems. 

We consider the following geometric situation.  Let $\pi\: Y \to S$ be a complex line bundle or a disk bundle over a
smooth base $S$ and $R \subset  Y$ a subspace such that $\pi \: R \to S$ is a covering space. Consider the category 
$\op_{R}(Y) = \{A \in \op(Y)|\text{Char}(A) \subset T^{*}_{Y}Y \cup T^{*}_{R}Y\}$ and the class of objects
${\Cal N} = \{A \in \op(Y)|\text{Char}(A) \subset T^{*}_{Y}Y\}$, i.e., the class of perverse sheaves that are local
systems on $Y$.  We form the localized category
$\oc_{R}(Y) = \op_{R}(Y)_{\Cal N}$. There exists a natural projection functor
$\Pi\:\op_{R}(Y)\to \oc_{R}(Y)  $.

\proclaim{Proposition 3.1} There exists a fully faithful functor $\oc_R(Y)\to \op_R(Y)$ which is a left quasi-inverse to
the projection
$\Pi\: \op_{R}(Y)\to \oc_{R}(Y) $. In the other words, there exists an embedding $\oc_{R}(Y) \to \op_{R}(Y)$ as a full
subcategory subcategory by a section of $\Pi$. \endproclaim

\demo{Proof} Let $\oc'_R(Y)$ be the full subcategory of $\op_R(Y)$ consisting of those $A$ for which
$R\pi_*(A)=0$. We claim that the restriction of $\Pi$ to $\oc'_R(Y)$ is an equivalence of categories $\Pi\:\oc'_R(Y)\to
\oc_R(Y)$. This is proved in Proposition 3.1 of
\cite{GMV}.
\enddemo

\proclaim{Proposition 3.2}  The functors $U \longmapsto \oc_{R}(\pi^{-1}(U))$ define a  stack on $S$. 
\endproclaim

\demo{Proof} Using the embedding $\oc_{R}(U) \hookrightarrow \op_{R}(U)$ from Proposition 3.1 we can identify
$\oc_{R}(U)$ with the full subcategory of $\op_{R}(U)$ which is characterized by the property that $R\pi_{*}A =
0$.  We see easily that since 
$\op_{R}$ is a stack the subcategories, $\{A \in \op_{R}(U)|R\pi_{*}A = 0\}$ form a stack on $S$.  Therefore,
$\oc_{R}$ is a stack.
\enddemo

We will now extend the previous considerations to the case when $R$ and the map $R \to S$ is allowed to have mild
singularities.

Let ${\overline S}$ be a smooth algebraic variety, $S\subset{\overline S}$ open dense subset, and assume that
${\overline S}-S$ is smooth and of codimension one.  Let
${\overline \pi}\: {\overline Y} \to {\overline S}$ be a line bundle and $\pi\:Y\to S$ its restriction to $S$. Let
${\overline R}
\hookrightarrow {\overline Y}$ be a divisor such that $\pi|R\:R \to S$ is \'{e}tale, $R = Y\cap {\overline R}$ and
${\overline
\pi}|{\overline R}\: {\overline R} \to{\overline S}$ is finite. We assume further that
${\overline R}$ has a Whitney stratification ${\overline R} = R \cup \bigcup^{k}_{i=1} {\overline R}_{i}$, where
all the strata
${\overline R}_{i}$ are of codimension one and the maps ${\overline \pi}|{\overline R}_{i}\: {\overline R}_{i} \to
{\overline S}-S$ are
\'{e}tale. Denote $\overline S_i =\overline\pi (\overline R_i)$.

Introduce the category $\op({\overline Y})$ of perverse sheaves $A$ such that
$\text{Char}(A) \subset T^{*}_{\overline Y} {\overline Y} \cup
\overline{T_{R}Y} \cup T^{*}_{{\overline R}_{i}}{\overline Y}$ and let ${\overline {\Cal N}}$ be the category
of local systems on
${\overline Y}$.  Denote $\oc_{\overline R}({\overline Y}) = \op({\overline Y})_{\overline {\Cal N}}$ the localized
category. Let also
$\oc_{R}({\overline Y})$ be the full subcategory of $\oc_{\overline R}({\overline Y})$ consisting of sheaves $A$
such that
$\text{Char}(A) \subset T^{*}_{Y}{\overline Y}
\cup \overline{T^{*}_{R}Y}$.

For $j:\ Y \hookrightarrow {\overline Y}$ we have the restriction functor $j^{*}:\ \oc_{\overline R}({\overline Y})
\to \oc_{R}(Y)$. If
$A \in \oc_{R}({\overline Y})$, then $A \cong j_{!*}j^{*}A$.  Therefore we can consider
$\oc_{R}({\overline Y})$ as a full subcategory of $\oc_{R}(Y)$.

To characterize the subcategory  $\oc_{R}({\overline Y})$ of $\oc_{R}(Y)$ we  proceed as follows. Choose a point
$y \in \overline R_i$ and let $x = \overline \pi (y) \in \overline S_i$. Choose a small neighborhood $U$ of $x$ and a
small neighborhood $V=U\times D'$ such that $V\cap \overline \pi^{-1}(x) \cap \overline R = \{y\}$.  Choose
$\epsilon \in D'$ such that $U \times \{\epsilon\} \cap
\overline R = 
\emptyset$. Let  $ j:\ (U-{\overline S}_{i})\times (D'-\{\epsilon\}) \hookrightarrow (U- {\overline S}_{i})\times D'.
$  For any $A
\in \op_{R}(Y)$ we can now consider the sheaf
$\var_{U,y}(A) = R^{-d}\pi_{*}j_{!}j^{*}A,\ d = \dim_{\Bbb C}S$.  It is a local system on $U - {\overline
S}_{i}$ and
$R^{-k}\pi_{*}j_{!}j^{*}A = 0$ for $k \neq d$. We call $\var_{U,y}(A)$ the {\sl variation} of $A$ on $U$ at the
point $y$. If $A$ is a local system then $\var_{U,y}(A) = 0$ and therefore $\var_{U,y}$  is defined on
$\oc_R(Y)$. We say that $A\in \oc_R(Y)$ has no variation at $y$ if  the the sheaf
$\var_{U,y}(A)$ is constant on $U-\overline S_i$ for small $U$. It is evident that it suffices to verify the absence of
variation on small $U$ such that $U$ is a standard ball in which 
$\overline S_i$ is given by the equation $x_d = 0$.

\proclaim{Proposition 3.3} The subcategory $\oc_R(\overline Y)$ of $\oc_R(Y)$ consists of sheaves $A \in \oc_{ R}(
Y)$ that have no variation at any point $y\in\overline R-R$.  
\endproclaim

\demo{Proof} Denoting by $k:Y \hookrightarrow \overline Y$ the inclusion  we have, as remarked earlier, that $A
\in \oc_R(\overline Y)$ necessarily implies $A \cong k_{!*}k^*A$. Let $A \in \oc_R(\overline Y)$. Consider the
following diagram of spaces and maps where
$y$ and $U$ have been chosen as above. 
$$
 \alignat 3 (U-{\overline S}_{i})\times (D'-t) \overset j\to\hookrightarrow & \ &(U-{\overline S}_{i})\times D'
&\overset {\overline j}\to\hookrightarrow  &\ &U\times D' \overset {\overline{\overline j}}\to\hookleftarrow U
\times (D'-\{\epsilon\})\\  &\ &\downarrow
\pi &\ &\ &\downarrow \pi \\  &\ &U-{\overline S}_{i} &\hookrightarrow &\ &U 
\endalignat 
$$ Since $A \in C_{R}({\overline Y})$ we have that $\text{Char}({\overline {\overline
j}}_{!}{\overline{\overline j}}^{*}A) \subset
\overline{T^{*}_{R}Y}
\cup\overline{T^{*}_{U\times\{\epsilon\}}{\overline Y}}\cup T^{*}_{\overline Y} {\overline Y}$, and therefore
$\text{Char}(R{\overline \pi}_{*}{\overline {\overline j}}_{!}{\overline{\overline j}}^{*}A)\subset
T^{*}_{U}U$, which means precisely that
$\text{var}_{U,y}(A) = 0$.

It remains to show that any $A\in \oc_R(Y)$ that has no variation lies in $\oc_{R}({\overline Y})$. Given such an
object $A$, we choose a  representative  $A$ in $\op_{R}(Y)$  such that $R\pi_{*}A = 0$.

Denoting  $k:\ Y \hookrightarrow {\overline Y}$, we claim that 
$\text{Char}(k_{!*}A) \subset T^{*}_{\overline Y}{\overline Y} \cup
\overline{T^{*}_{R}(Y)}$. We must show that the components $T^{*}_{{\overline R}_{i}}{\overline Y}$ and
$T^{*}_{{\overline
\pi}^{-1}({\overline S}_{i})}{\overline Y}$ do not appear in $\text{Char}(k_{!*}A)$.  This is clearly a local
question and by cutting with a normal slice to ${\overline S}_{i}$ we can assume that
$U \subset {\Bbb C}$ and ${\overline S}_{i} = \{0\}$.

The assumption $R\pi_{*}A = 0$, together with the fact that
$R^{k}\pi_{*}j_{!}j^{*}A = 0$ for $k \neq -1$ $(1 = \dim S)$, implies that $A|U \times
\{\epsilon\}$ is constant.  This implies that the components of the form
$T^{*}_{{\overline \pi}^{-1}({\overline S}_{i})}{\overline Y}$ do not appear in $\text{Char}(j_{!*}A)$.

The fact that the components $T^{*}_{{\overline R}_{i}}({\overline Y})$ do  not appear follows now directly
from the fact that the variation map is zero (see, e.g., \cite{MV1, Theorem 5.3}). 
\enddemo

\heading 4. Melded Morse systems
\endheading
                                
In this section we construct a stack $\om$ on $\ct$ and a morphism of stacks
$\ode\to\om$. Just as objects of $\ode$ have Lagrangian support, so do the objects of $\om$ and they are indeed
defined by first fixing this support to be $\Lambda$, a given conic  Lagrangian subvariety of $T^*X$. In Section 5
we prove  that this morphism is an  equivalence outside a codimension 2 locus in the support we have fixed. In
section 6 we show that  this morphism is an embedding, i.e., for any $U\subset \ct$, $\ode(U)
\to \om(U)$  is a fully faithful functor.

\subheading{Families of Morse Functions}

We want to produce families of Morse functions (or test triples).  Let us consider $T\ct$  as a bundle over $\ct$, and
let $Gr(T\ct)$ be the associated bundle of Lagrangian Grassmannians. Let $Gr^0(T\ct)$ be the open part of
$Gr(T\ct)$ consisting of {\sl horizontal}  planes, i.e. planes on which $d\pi\:TT^*X\to T^*X$ is an isomorphism. 
(Planes in  $Gr^0(T\ct)$ are transverse to the  tangent plane to the fiber of $\pi$.) If $\rho$ denotes the projection
$Gr^0(T\ct)\to X$, then $\rho^{-1}(x)$ is canonically identified with the space of $2$-jets of complex valued
functions, modulo constant functions, at $x$. Let $B_\delta(x) \subset X$ denote the ball of radius $\delta$ around
$x\in X$, as  measured by some fixed Riemannian metric, and $B_\delta^0(x)$ denote its interior.

\proclaim{\bf Lemma 4.1} To each $L\in Gr^0(T\ct)$, we can associate a
$\delta(L)\in\Bbb R$ and an analytic function
$\phi_L:B_{\delta(L)}^0(\rho(L))\to \Bbb C$ in such a way that $L$ is  the tangent space to the Lagrangian
manifold $\{d\phi_L\}$.  Furthermore, we can make these choices in such a way that they vary  smoothly with
respect to $L$ in the following sense:  Let $W$ be the union the disks
$B_\delta^0(\rho(L))\cross L$ in $X\cross Gr^0(T\ct)$, and let $f$ be the  function $(x,L)
\mapsto \phi_L(x) $.  Then $W$ is open in $X$ and is bounded by a smooth  (not analytic!) submanifold, and $f$ is
a smooth function on $W$.   Furthermore, $f$ can be chosen to be invariant under the circle group $S^1\subset \Bbb
C^*$, where $S^1$ acts on
$Gr^0(T\ct)$ through its action on $\ct$, and acts on functions by multiplication of the values. \endproclaim

\demo {Proof}  To do this when $X$ is Euclidean space is easy:  $\phi_L$ can  be taken to be the unique polynomial
of degree $2$ such that
$\phi_L(\rho(L))=0$ and $L$ is the tangent space to the Lagrangian manifold
$\{d\phi_L\}$.  Now the general case is obtained from a covering by  coordinate charts by a smooth partition of
unity argument.  
\enddemo

\subheading{The construction of the stack $\om$}   Let $\Lambda \subset X$ be a  conic closed complex
Lagrangian subvariety.  We write $\ode_\Lambda$ for the full substack of $\ode$ consisting of objects whose
characteristic variety lies in $\Lambda$.   We will define the $\om_\Lambda$ of melded Morse systems.   The stack
$E$ is $\ode = \bigcup
\ode_\Lambda$ where the union is taken over all Lagrangian subvarieties $\Lambda$; an analogously for the stack
$\om$.    

The construction of $ \om_\Lambda (U)$ proceeds in two steps, corresponding to the codimension zero strata and
the codimension one strata of $\Lambda$. 
\medskip

\noindent{\bf Step 0.}   We define $\Lambda^0$ to be the codimension $0$ stratum of $\Lambda$.   Specifically, let
$$
\align
\Lambda^0 = \{\xi\in\Lambda | \xi \text { a smooth point and }& T_{\xi'}\Lambda \cap T_{\xi'}\pi^{-1}(\pi(\xi'))
\text{ is of constant}\\ &\text{dimension for } \xi'\in\Lambda \text { near } \xi \}. 
\endalign
$$ 

The restriction of the stack $\om$ to $\Lambda^)$ is just the stack of local systems on $\Lambda^0$.  (This is just
what you would expect, since by the vanishing cycle construction, a Perverse sheaf constructible with respect to
$\Lambda$ give rise to a local system on $\Lambda^0$.)  However, we will define $\om$ on $\Lambda^0$ by an
equivalent definition that is more complicated.  The reason for this is to facilitate the extension of the stack over the
rest of $\Lambda$, as described in Step 1 below.

We define a fiber bundle $\tau^0 : Q^0 \to \Lambda^0$ over $\Lambda^0$ as follows: 
$$ Q^0 = \{ L \subset T_\xi (T^*X) | L \text{ Lagrangian }, L \pitchfork  T_\xi\pi^{-1} (\pi(\xi)) , L \pitchfork
T_\xi\Lambda \}.
$$ We have a canonical map 

$$ Q^0 \to Gr(T\ct) \times Gr(T\ct) \times Gr(T\ct)
$$ given by sending the point $(\xi, L\subset T_\xi\ct) \in Q^0$ to
$(L,T_\xi\pi^{-1}(\pi(\xi)), T_\xi\Lambda)$. Pulling back the Maslov bundle on the appropriate subset of $Gr(T\ct)
\times Gr(T\ct)
\times  Gr(T\ct)$ we get the Maslov bundle $\lambda$ on $Q^0$.

We define two families of perverse sheaves on a disk, both constructible with respect to the zero of the trivial
bundle, as described in \S 3.  The first is parameterized by $\Lambda^0$

$$ \oc_{\Lambda^0 \times \{ 0 \} } (\Lambda^0 \times \Bbb C)$$

and the second is parameterized by $Q^0$

$$ \oc_{Q^0 \times \{ 0 \} } (Q^0 \times \Bbb C)$$

The circle group $S^1$ acts on $Q^0 \times \Bbb C$ as follows:  it acts on $\ct$ by scalar multiplication
of vectors and it actos on the complex numbers $\Bbb C$ by scalar multiplication.

We define the stack $\om^0_\Lambda$ on $\ct$ to be

$$ \om^0_\Lambda = \tau^0_*(\oc_{Q^0 \times \{ 0 \} } (Q^0 \times \Bbb C)^{S^1})_\lambda$$

where the superscript $S^1$ denotes the full substack of $S^1$ equivariant objects and the subscript $\lambda$
denotes the full substack consisting of objects with the property that if we tensor them with the local system
$\lambda$ they descend to objects of
$\oc_{\Lambda^0 \times \{ 0 \} } (\Lambda^0 \times \Bbb C)^{S^1}$.

\medskip
\noindent{\bf Remark.}  The stack $\om^0_\Lambda$ has three other equivalent descriptions

$$ \om^0_\Lambda = \oc_{\Lambda^0 \times \{ 0 \} } (\Lambda^0 \times \Bbb C)^{S^1} = \{ \text{Local systems
on} \ Q^0 \}_\lambda = \{
\text{Local systems on} \ \Lambda^0 \} $$

Next, we choose a family of test manifolds ${\Cal T }_q\in T^*X $  parameterized by points $q\in Q^0$ such that if
$q$ represents a subspace $L\subset T_\xi  T^*X$, then
${\Cal T }_q$ is tangent to $L$. To be precise, for each $L\in Q^0$ we choose a triple
$(B_{\delta(L)}(\rho(L)),\phi_L, D_{\epsilon(L)})$ as follows. The pair
$(B_{\delta(L)}(\rho(L)),\phi_L)$ is chosen as in the lemma on families of  Morse functions above, so $\phi_L
(\rho(L))=0$, and (this is the main point) $L$ is  tangent to
$\{d\phi_L\}$. The triple must be a test triple. We further require that for all $\delta<
\delta(L)$ there exists an $\epsilon$ such that $(B_{\delta}(\rho(L)),\phi, D_{\epsilon})$ is still a test triple. These
conditions can be met for each $L$ by Lemma 3.5.1 from
\cite{GM2}; in fact, the set of pairs $\delta(L), \epsilon(L)$ satisfying the conditions is a fringed set $0<\epsilon\ll
\delta \ll 1$ \cite{GM2, \S 5}. That we can choose
$\delta(L),\epsilon(L)$ globally and smoothly follows from \cite{GM2, \S 5.5}.

To define the morphism of stacks $\nu^0\: \ode_\Lambda  \to \om_\Lambda^0$, we first note that it suffices to
define the functor 
$\nu^0(U)\: \ode_\Lambda(U)  \to
\om_\Lambda^0(U)$ on small  open sets $U \subset \ct$. By the construction of $\ode$, it suffices to define the
functor $\nu^0(U)$ on
$\op(X,U)$. This functor is defined on objects  as follows. Let $W$ be the union of the disks
$B_{\delta(L)}^0(\rho(L)) \intersect
\phi^{- 1}(D^0_{\epsilon(L)}  \cross \{L\})$ in $X\cross Q^0$, and let $f:W\to Q^0\cross \Bbb C$ be the  function
$(x,L) \mapsto (L,\phi_L(x)) $. Let $pr$ be the projection of $W$ to $X$. For $A\in\op_\Lambda (X,U)$ we 
consider
$R^{-1} f_*pr^*A$ as an object in $\oc_{Q^0\times\{0\}}(Q^0\times\Bbb C)$. 

\proclaim{Proposition 4.2} The construction above is independent of the  choices and it defines a functor
$\nu^0(U)\: \op_\Lambda (X,U) \to  \om_\Lambda^0(U)$. 
\endproclaim 

Note that in defining the functor $\nu^0(U)$ we first pick  $A\in\op_\Lambda  (X,U)$ and then we choose $W$
accordingly. 

\demo {Proof}The group $S^1$ acts on $Q^0 \times \Bbb C$, $W \subset X \times Q^0$, and $X$, where the action
on $X$ is trivial.  The maps $pr$ and $f$ are equivariant with respect to this action.  Hence $R^{-1} f_*pr^*A$ lies
in  $\oc_{Q^0\times\{0\}}(Q^0\times\Bbb C)^{S^1}.$

We have 

$$ \{ \text{Local sytems on} \ Q^0\} = \oc_{Q^0 \times \{0\}} (Q^0 \times \Bbb C)^{S^1}$$

Let $\Cal L$ be the local system on $Q^0$ corresponding to $R^{-1} f_*pr^*A$.  The fact that ${\Cal
L}\otimes\lambda$ descends to a local system on $\Lambda^0
\intersect U$ may be seen as follows. The stalk of ${\Cal L}$ is the local Morse group of \cite{GM2, p.224}. Since
local Morse data is tangential Morse data tensored with normal Morse data, we have that ${\Cal L}={\Cal
L}_T\otimes {\Cal L}_N$ by the K\"unneth theorem. (Since we are dealing with perverse sheaves, there is only one
nonvanishing cohomology group in each case.) By the classical Morse theory, ${\Cal L}_T$ is one-dimensional and
coincides  with the  Maslov line bundle $\lambda$. Since the Maslov line bundle is isomorphic to its inverse,
$\lambda=\lambda^{-1}$, we have
${\Cal L}\otimes\lambda={\Cal L}_N$. But ${\Cal L}_N$ is the vanishing cycle cohomology, which lives on
$T^*X$.   

Given two sets of choices, we use the same technique of \cite{GM2,\S 5.5} to find an one-parameter family of
choices connecting them.  This shows  that the resulting local system ${\Cal L}$ is independent of the choices.   
\enddemo

\medskip \noindent {\bf Step 1.} We define $\Lambda^1$ to be the codimension $1$ stratum of
$\Lambda$. Specifically, let $\Lambda^1 $ be the largest subset of $\Lambda-
\Lambda^0$ such that $T_{\xi'}\Lambda^1 \cap T_{\xi'}\pi^{-1} (\pi(\xi'))$ is of constant dimension for $\xi'$ near
$\xi$ in
$\Lambda^1$ and
$\Lambda^0\union\Lambda^1$ is Whitney stratified.  

We define a space $\tau^1\:Q^1\to\Lambda^1$ over $\Lambda^1$ as follows:  $Q^1$ is the space of triples
$(\xi,L,\eta)$, where
$\xi\in\Lambda^1$; $L$ is a Lagrangian subspace of $T_\xi T^*$X such that $L\intersect T_\xi\Lambda^1 =\{0\}$,
$L\pitchfork T_\xi\pi^{-1} (\pi(\xi))$ and $L$ is transverse to all Lagrangian subspaces $L'\subset T_\xi T^*X $
that are limits of tangent  spaces to points in $\Lambda^0$; and $\eta \in T_\xi T^*X/(L\oplus  T_\xi\Lambda^1)$ is
a nonzero vector. (Note that $ T_\xi T^*X/(L\oplus  T_\xi\Lambda^1)$ is one-dimensional.)

Note that by making $\Lambda^1$ smaller, we could achieve that $\tau^1\:Q^1 \to
\Lambda^1$ would be a fiber bundle; however this offers an aesthetic advantage only.

Next, we choose a family of test manifolds ${\Cal T}_q\in T^*X$  parameterized by points $q\in Q^1$ such that if
$q=(\xi,L,\eta)$, then ${\Cal T}_q$ is  obtained by displacing a test manifold tangent to $L$ by a small
displacement in the  direction of
$\eta$.

We let $\overline{Q}^1 $ denote the space of triples  $(\xi,L,\eta) $ as in $Q^1$, but  allowing $\eta$ to be zero. 
First, we note that the line $T_\xi T^*X/(L\oplus  T_\xi\Lambda^1)$ is canonically identified with a complement to
$T_\xi\pi^{-1}(\pi(\xi))
\intersect L\oplus T_\xi\Lambda^1$  in  $T_\xi\pi^{-1} (\pi(\xi))$.By making a smoothly varying  choice of such
complementary lines, we can think of $\eta$ as lying in 
$T_\xi\pi^{-1} (\pi(\xi))$, which is canonically identified with $T^*_{\pi(\xi)}X$. Now for each  $q=(\xi,L,\eta)\in
\overline{Q} ^1$, we choose a plane $L(q,\kappa)\in Gr^0(T_{\xi+\kappa\eta} \ct)$, where $\kappa\in\Bbb R$ and
$\xi+\kappa\eta$ refers to addition in $T^* _{\pi(\xi)}X$. We make this choice in any continuous way, requiring
only  that if $\eta=0$ then $L(q,\kappa)=L$. (Note that
$\rho(L(q,\kappa))=\rho(L)$.)

For each $q=(\xi,L,\eta)\in Q^1$ we choose a triple  $(B_{\delta(L)}(\rho(L)),
\phi_{q,\kappa(L,\eta)},D_{\epsilon(L)})$, subject to the following conditions:
$\{d\phi_{(q,\kappa)}\}$ is tangent to $L(q,\kappa)$ and the triple  is a test triple. Furthermore, for all
$\delta<\delta(L)$, there should exist an $\epsilon$ such that
$(B_{\delta}(\rho(L)),\phi, D_{\epsilon})$ is still a test triple. 

The proof that this can be achieved is similar to the proof for $Q^0$ except  that now we have a fringed set of triples
$\delta(L),\epsilon(L),
\kappa(L,\eta)$ satisfying the conditions $0<\kappa \ll\epsilon\ll\delta\ll 1$. (The point is that if we have a test triple
$(B,f,D)$ and we perturb $f$ by a small amount $\kappa$, we still have a test triple.)

Let $W^1$ be the union the disks $B_{\delta(L)}^0(\rho(L))\intersect
\phi^{-1}_{q,\kappa}(D^0_{\epsilon(L)}\cross \{q\})$ in $X\cross Q^1$,  and let
$f:W^1\to Q^1\cross D$ be the function $(x,q) \mapsto (L,\epsilon^{-1} 
\phi_{q,\kappa}(x))$. Here $D\subset \Bbb C$ is the unit disk.  For  each point $q\in Q^1$, the singular values of
$\epsilon^{- 1}\phi_{q,\kappa}$ form a finite set of points, which varies smoothly with  $q$. Let us denote the set
of these values by $R$. Then
$R\subset Q^1 \times D$ and the projection $R\to Q^1$ is a topological covering. To define $\om_\Lambda$, we
first construct $\tilde Q^1\subset Q^0 \times Q^1$ as follows,
$$
\align
\tilde Q ^1 = \{(q^0,q)\in Q^0 \times Q^1 \mid{}& q^0 = (\xi^0,L^0), L^0 \subset T_{\xi^0}T^*X\\  &\text{such
that } \xi^0 \in \Cal T_q \cap \Lambda^0 \text{ and } T_{\xi^0}\Cal T_q = L^0 \}.
\endalign 
$$ Then $\tilde Q^1 \to Q^1$ is a topological covering canonically isomorphic to $R \to Q^1$. In particular we get a
map $r:\tilde Q^1 \times D \to Q^1  \times D$ such that
$\tilde Q^1 \times \{0\}$ maps isomorphically onto $R$ and  $\tilde Q^1 \times D$ maps to a small tubular
neighborhood of $R$ in $Q^1
\times  D$. We also have a natural map $s: \tilde Q^1 \times D \to Q^0 \times D$ which  is the projection $\tilde Q^1
\to Q^0$ on the first component and identity on  the $D$-component.

If we are given $A \in C_{Q^0 \times \{0\}}(Q^0 \times D)$, then $s^*A \in C_{\tilde Q^1 \times \{0\}}(\tilde Q^1
\times D)$ is well defined and if $B \in C_R(Q^1 \times D)$, then $r^*B \in C_{\tilde Q^1 \times \{0\}} (\tilde Q^1
\times D)$ is also well defined. We also note that we have an inclusion $Q^1 \subset \overline Q^1$ and we also
have $R \subset \overline R$, where $\overline R$ is defined, just like $R$ is, as the locus of the critical values of
$\epsilon^{-1}\phi_{q,\kappa}$ except that now we also allow $\eta = 0$. We are now in the situation of Section 3
and the category $ C_R(\overline Q^1 \times D)$ makes sense with respect to  $(Q^1,R)
\subset (\overline Q^1,\overline R)$. We are now ready to define the stack $\om_\Lambda$

\proclaim {\bf Definition 4.3} Let $U$ be an open subset of $\ct$.  Define  $ Q^0(U)$ to be $(\tau^0)^{-1} (U)
\subset Q^0$,
$\overline Q^1(U)$ to be $(\tau^1)^{-1} (U) \subset \overline Q^1$ and $\tilde Q^1(U)$ to be $r^{-1} \overline
Q^1(U) \cap s^{-1} Q^0 (U)$.  The category  $\om_\Lambda (U)$ is defined as follows. Its objects are triples
$(A,B,\gamma)$, where

\roster
\item  $A \in
\om^0_\Lambda(U) = C_{Q^0 (U) \times \{ 0 \}} ( Q^0(U) \times D)^{S^1},$ 
\item $B \in  C_{R(U)} (\overline Q^1 (U) \times D),$
\item $\gamma: s^*A @>\cong>> r^*B$ is an isomorphism in the category $C_{\tilde Q^1 (U) \times \{0\}}(\tilde
Q^1(U) \times D)$ 
\endroster The morphisms are defined in the obvious way.
\endproclaim

It is clear from this definition that $\om_\Lambda$ is a stack.

To define the morphism of stacks $\nu^1\: \ode_\Lambda \to \om_\Lambda$ we  proceed as follows. Again it
suffices to define the functor $\nu^1(U)$ for small open sets and we fix $P\in\op_\Lambda(X,U)$, and we make the
choices for this  particular $P$. Again as for $\nu^0$ all good choices give the same functor. Recall that $W^1
\hookrightarrow Q^1 \times X$  and write $pr\: W^1 \to X$ for the projection. Then for $P \in P_\Lambda(X,U)$
we  define $$
\nu^1(U)(P) = (\nu^0(U)(P),R^{-1}f_*pr^*P) \in C_{Q^0 \times \{0\}}(Q^0  \times D)^{S^1}\times  C_R(\overline Q^1
\times D).
$$ This evidently gives us morphism of stacks  $\nu_\Lambda\: \ode_\Lambda \to 
\om_\Lambda$, and passing to the limit over $\Lambda$ we get a morphism of stacks
$\nu\: \ode \to \om$.
      
\heading 5. Equivalence through codimension 1
\endheading

In Section 4 we constructed a stack $\om$ on $\ct$ and a morphism of stacks  $\nu\:
\ode \to \om$. The objects of both $\ode$ and $\om$ are supported on Lagrangian  subvarieties of $\ct$. Given 
$\Bbb C^*$-conic open subset $U \subset \ct$ and a  closed, $\Bbb C^*$-conic Lagrangian subvariety $\Lambda 
\subset U$, let
$\Lambda^0,\Lambda^1  \subset \Lambda$ be the subsets defined in Section 4, and let $\Lambda^2 = \Lambda 
-(\Lambda^0 \cup
\Lambda^1)$.

\proclaim{\bf Theorem 5.1} The functor $\nu\: \ode_\Lambda(U- \Lambda^2) \to
\om_\Lambda^1 (U-\Lambda^2)$ is an equivalence of categories. 
\endproclaim
                             
\demo{Proof} Because both $\ode$ and $\om$ are stacks it suffices to prove the theorem on stalks. Therefore we are
reduced to two cases: We consider stalks at $\xi \in \Lambda^0$ or $\xi \in
\Lambda^1$. We start with the case $\xi \in \Lambda^0$.

Let $\xi \in \Lambda^0$ and we choose a small open, $\Bbb C^*$-conic neighborhood
$U$ of $\xi$ such that $U\cap\Lambda^0/\Bbb C^*$ is isomorphic to a ball. Choose also
 a contact transformation $\phi : U \to V$, $V \subset \cty$ such that $\phi(\Lambda\cap U)$ is in generic position.
Without loss of generality we can assume that $Y = \Bbb C^n$, $\phi(\xi)=(0,dy_n)$ and that $V$ projects onto a
small open neighborhood $W$ of the origin. Consider the projection $\Bbb C^n \to \Bbb C^{n-1}$ given by $(y_1,
\dots, y_n) \mapsto (y_1,\dots,y_{n-1})$. Shrink $W$ so that it is of the form
$W=S\times D$, where $S\subset \Bbb C^{n-1}$ and $D$ is a disk in $\Bbb C$. Then we get a projection $W \to S$ with fibre $D$. We see
 by Lemma 2.5 that the category $\ode_\Lambda(U)$ is equivalent to $\oc_{S\times\{0\}}(S\times D)$, and therefore
it is equivalent to the category of local systems of finite rank on $U$. 

On the other hand, if we fix $s\in S$ then via the contact  transformation $\phi$ the set $\{s\} \times D$
parameterizes a family of hypersurfaces in $X$ and therefore it gives a map from an open set $B \subset X$ to $D$.
In other words, it gives us a melding. Altogether, then, $S$ parameterizes a family of meldings giving us a map $S \to
Q^0$. The category of meldings also gives us the category of local systems of finite rank on $U$ and via our
identification of $\ode_\Lambda(U)$ with $\oc_{S\times\{0\}}(S\times D)$ we see that the functor $\nu^0
:\ode_\Lambda(U) \to
\om^0_\Lambda(U)$ becomes essentially identity.

Let us now assume that $\xi \in \Lambda^1$. We again choose a small open, $\Bbb C^*$-conic neighborhood $U$
of $\xi$ and a contact  transformation $\phi : U \to V$, $V \subset \cty$ such that
$\phi(\Lambda\cap U)$ is in generic position.  Without loss of generality we can assume that $Y = \Bbb C^n$,
$\phi(\xi)=(0,dy_n)$ and that $V$ projects onto a small open neighborhood $W$ of the origin. The situation is now
slightly more complicated than in the case of $\Lambda^0$. The components of $\phi(\Lambda)$ project to
hypersurfaces $R_i$ in $\Bbb C^n$ which all have $y_n = 0$ as the tangent cone at the origin. Because the
components of $\Lambda$ are tangent to each other along $\Lambda^1$ we see that $\Lambda^1$ projects to a
codimension 2 submanifold, which we can assume is given by $\Delta = \{(y_1, \dots,y_n) | y_n = y_{n-1} = 0 \}$,
and all the hypersurfaces $R_i$ are tangent along $\Delta$. Let us write $\bar R = \bigcup R_i$ and $R = \bar R -
\Delta$. Consider the projection $\Bbb C^n \to \Bbb C^{n-1}$ given by $(y_1, \dots, y_n) \mapsto
(y_1,\dots,y_{n-1})$. Shrink $W$ so that it is of the form
$W=\bar S\times D$, where $\bar S\subset \Bbb C^{n-1}$ is an open ball around the origin and $D$ is a disk in
$\Bbb C$.  Then we get a projection $W
\to \bar S$ with fibre $D$. Let $S = \bar S -
\{y_{n-1}=0\}$. Then, by shrinking $\bar S$ if necessary, we get into the situation of section 3. We have a
projection $\pi: \bar S
\times D \to \bar S$ and
$\bar R \subset \bar S$ such that
$\pi^{-1}(S) \cap \bar R = R$ and $ \pi : R \to S$ is a topological covering. 

Note that in the process above we have been shrinking the open set $V \subset \cty$. Let us denote again by
$U\subset \ct$ the set corresponding to the modified $V$ via $\phi$. Then Lemma 2.5 tells us that
$\ode_\Lambda(U)$ is equivalent to $\oc_R(\bar S \times D)$. Recall that in section 3 we gave an explicit
construction of the category $\oc_R(\bar S \times D)$ as a subcategory of  $\oc_R(S
\times D)$. We can view $S\times D \to S$ as a family of meldings and as such we get a map $S \to Q^1$. In this
way we get a functor
$\tau :
\om_\Lambda (U) \to \ode_\Lambda(U)$ by pulling back the melding data in part (2) of Definition 4.3 from $Q^1$
to $S$. It is clear that $\tau \circ
\nu$ is equivalent to identity. Therefore, it remains to show that $\nu \circ \tau$ is naturally equivalent to identity.

To this end let us consider $(A,B,\gamma) \in \om_\Lambda (U)$, i.e.,  $A \in  \oc_{Q^0 \times \{0\}}(Q^0 \times
D)^{S^1}$, $B \in  C_{R^1}(\bar Q^1 \times D)$, and $\gamma: s^*A @>\cong>> r^*B$. We have here denoted by
$R^1$ the set which was denoted by $R$ in definition 4.3. Let us denote by $B'$ the image of $B$ under the
restriction functor $\oc_{R^1}(\bar Q^1 \times D) \to \oc_R(\bar S
\times D)$. We must now show that the triple
 $(A,B,\gamma)$ can be reconstructed from $B'$. Let us denote $W = \bar S \times D$ as we did above. Then the
vanishing cycle functor
$$ \oc_R(\bar S \times D) \to
\{\text{local system on } T^*_RW-T^*_WW \} $$ gives us $A$ by the first part of the proof. The fact that
$B'$ also determines $B$ and $\gamma$ follows from Lemma 5.3 below which we now set up.

\enddemo
\definition{Definition 5.2}  Let $\Cal C_0$, $\Cal C_1$, $\Cal C_2$ be three categories,  $\alpha\:\Cal C_1\to 
\Cal C_0$, $\beta\:\Cal C_2\to  \Cal C_0$ two functors. The {\it fiber product\/} $\Cal C=\Cal C_1\times _{\Cal
C_0} \Cal C_2$ is defined as the category defined as follows. Objects of $\Cal C$ are triples $(A,B,r)$, where $A\in
\Ob \Cal C_1$, $B\in\Ob \Cal C_2$,
$r : \alpha(A) \cong \beta(B)$ is an isomorphism in $\Cal C_0$. Morphisms $(A,B,r)\to (A',B',r')$ in $\Cal C$ are
pairs of morphisms
$A\to A'$, $B\to B'$ satisfying a natural commutativity conditions with $r$ and
$r'$.
\enddefinition

Mappings $(A,B,r)\mapsto A$ and  $(A,B,r)\mapsto B$ give rise to functors
$\overline\alpha\Cal C\to \Cal C_2$ and $\overline\beta\Cal C\to \Cal C_1$, and $r$ defines the isomorphism of
compositions 
$\alpha\overline \beta \cong \beta \overline
\alpha$. Moreover, given a category $\Cal D$ and diagram of functors 
$$
\CD
\Cal D @>\widetilde\beta>> \Cal C_1\\ @V\widetilde\alpha VV     @VV\alpha V\\
\Cal C_2  @>\beta>>  \Cal C_0
\endCD
$$ together with an isomorphism of functors $\alpha\widetilde\beta \cong
\beta\widetilde\alpha$, there exists a unique functor $\delta\:\Cal D\to \Cal C$ such that
$\widetilde\alpha = \overline\alpha \delta$, $\widetilde\beta = \overline\beta \delta$.

We use this construction in the case
$$
\align
\Cal C_1&=\{\text{local systems on $T_R^*(S \times D)-T^*_{S\times D}(S\times D)$}\},\\ 
\Cal C_2&= C_{\pi^{-1}(s)}(\{s\}\times D),\\ 
\Cal C_0&=
\{\text{local systems on $T_{\pi^{-1}(s)}^*(\{s\} \times D) -  T^*_{\{s\}\times D}(\{s\}\times D)$}\},
\endalign
$$ 
$$
\align
\alpha\:\Cal C_1\to \Cal C_0 &\quad \text{is the restriction of local systems,}\\
\beta\:\Cal C_2\to \Cal C_0 &\quad \text{is the vanishing cycle functor.} \endalign
$$

Take $\Cal D=\oc_R(S\times D)$. Then the functors
$$
\widetilde\alpha\:C_R(S\times D)\to  \oc_{\pi^{-1}(s)}(\{s\}\times D) \tag1
$$ (restriction) and
$$
\widetilde\beta\:\oc_R(S\times D)\to \{\text{local systems on $T_R^*(S \times D) - T^*_{S\times D}(S\times D)$}\}
\tag2
$$ (vanishing cycles) satisfy the condition that  $\alpha\widetilde\beta$ is isomorphic to
$\beta\widetilde\alpha$. Therefore, we get a functor $\delta\:\oc_R(S\times D)\to
\Cal C$. 

\proclaim{Lemma 5.3} Let $S$  be connected and $s \in S$. Then the functor  
$$ 
\delta\:\oc_R(S \times D) \to \Cal C 
$$   is fully faithful, i.e., $\oc_R(S \times D)$ is a full subcategory of $\Cal C$. \endproclaim

\demo{Proof} The proof is based on technique of the action of a group on a category used in \S4 of
\cite{GMV}.

Let us number the points of $\pi^{-1}(s)\cap R$,  $\pi^{-1}(s)\cap R=\{x_1,\dots,x_n\}$. Recall that in \cite{GMV}
we defined the category $\Cal Q$ of quivers $Q=(M_i,m_{ji})$, where $M_i$, $1\leq i\leq n$, are
finite-dimensional linear spaces, $m_{ji}\:M_j\to M_i$ are linear maps such that
$1+m_{ii}$ are invertible.

Denote by $G$ the fundamental group of $S$, $G=\pi_1(S,s)$. The group $G$ acts on the $n$-point set
$\pi^{-1}(s)\cap R$, and in
\cite{GMV} we defined the compatible action $\Phi$ of $G$ on the category $\Cal Q$ such that $\oc_R(S\times D)$
is equivalent to $\Cal Q_\Phi$, the category of $G$-equivariant objects of $\Cal Q$. Recall that objects of  $\Cal
Q^\Gamma$ are families $\{Q, \rho(g)\}$, where  $Q\in \Ob \Cal Q$, and
$\rho(g)\:Q\to \Phi(g) Q$, $g\in G$, are  isomorphisms in $\Cal Q$.

Now we present a similar quiver description of the category $\Cal C_1$ of local systems on $T_R^*(S \times
D)-T^*_{S\times D}(S\times D)$.

We have the projection map  $T_R^*(S \times D)-T^*_{S\times D}(S\times D)\to S$. Clearly, the fiber $X$ of $p$
over $s\in S$ is the space
$T_{\pi^{-1}(s)}^*(\{s\} \times D) -  T^*_{\{s\}\times D}(\{s\}\times D)$.

\proclaim{Lemma 5.4} $X$ is homeomorphic to the disjoint union of $n$ punctured disks. In particular, the
category $\Cal C_0$ of local systems on $X$ is equivalent to the category $\Cal L$ of quivers $L=(L_i,l_i)$, where
$L_i$, $i=1,\dots,n$, are finite dimensional linear spaces and $l_i\:L_i\to L_i$ are invertible linear maps.
\endproclaim

\demo{Proof} Clear, since
$$ X=\bigcup_{x_i\in \pi^{-1}(s)\cap R} \Bbb C^*.
$$ 
\enddemo

The action of $G$ on  $\pi^{-1}(s)\cap R=\{x_1,\dots,x_n\}$ defines the action $\Psi$ of $G$ on the category $\Cal
L$:
$$
\Psi(g)(L_i,l_i)=(L_{g(i)},l_{g(i)}).
$$

\proclaim{Lemma 5.5} The category $C_1$ of local systems on $T_R^*(S \times D)- T^*_{S\times D}(S\times
D)$ is equivalent to the category $\Cal L_\Psi$ of $G$- equivariant objects for the above action of $G$ on $\Cal L$.
The functor $\alpha\:\Cal C_1\to \Cal C_0$ is given by $(L,\rho_L(g))\to L$. 
\endproclaim

\demo{Proof} The proof is a simplified version of the proof of Theorem 4.7 in
\cite{GMV}.
\enddemo

Using this equivalence of categories, we represent the functor
$\widetilde\beta$ in (2) as follows.

Let we are given an object $\{Q,\rho(g)\}$ in $\Cal Q_\Phi$, so that we have an object
$$ Q=(M_i,m_{ji})
$$ in $\Cal Q$ and isomorphisms
$$
\rho(g)\:Q\to \Phi(g)Q
$$ in $\Cal Q$, represented by isomorphisms
$$
\gamma_{i,g}\:M_i \to M_{g(i)}
\tag3
$$ of linear spaces such that for each $i,j,g$ the diagram
$$
\CD M_i   @>\gamma_{i,g}>>  M_{g(i)}\\ @Vm_{ji}VV   @VV\widetilde m_{g(j),g(i)}V \\ M_j
@>\gamma_{j,g}>> M_{g(j)}
\endCD
\tag4
$$ commutes. Here $\widetilde m_{g(j),g(i)}$ are certain polynomials in $m_{kl}$ and
$(1+m_{kk})^{-1}$ constructed using the homomorphism of the fundamental group
$G=\pi_1(S,s)$ to the braid group $B_n$, as explained in \cite{GMV, Proposition 2.4}. In particular formulas
therein show that 
$$
\widetilde m_{g(i)g(i)}= m_{g(i)g(i)}.
\tag5
$$ 

Given an object $\{Q=(M_i,m_{ji}),\rho(g)\}$ in $\Cal Q_\Phi$ we define the object 
$\{L=(L_i,l_{ji}),\allowmathbreak 
\rho_L(g)\}$ in $\Cal L_\Psi$ setting $L_i=M_i$,
$l_i=1+m_{ii}$, and defining $\rho_L\:L\to \Psi(g)L$ by the same linear maps as the morphism $\rho(g)\:Q\to
\Phi(g)Q$. The diagrams (4) for $i=j$, together with formula (5), show that $\{L,\rho_L\}$ is indeed an object from
$\Cal L_\Psi$. In these notations, the functor $\widetilde\alpha$ from (1) is given by $(Q,\rho(g))\mapsto Q$,  
$\widetilde\beta$ from (2) is given by
$\{Q,\rho(g)\}\mapsto  \{L,\rho_L(g)\}$, and the functor 
$$
\delta\:\Cal Q_\Phi\to \Cal L_\Psi\times _{\Cal L} \Cal Q
$$ from Lemma 5.3 is given by
$$
\delta\:\{Q,\rho(g)\}\mapsto  (\{L,\rho_L(g)\}, Q, \operatorname{id}_L). $$

We want to prove that $\delta$ is fully faithful, i.e., that it induces isomorphisms of all
$\Hom$-spaces.

Clearly, $\delta$ is injective on $\Hom$-spaces.

To prove that $\delta$ is surjective on $\Hom$-spaces, take objects $\{Q,\rho(g)\}$,
$\{Q',\rho'(g)\}$ in $\Cal Q_\Phi$. Let $\tau$ be a morphism $$
\tau\:\delta(\{Q,\rho(g)\}) \to \delta(\{Q',\rho'(g)\}).
$$ We must construct $\sigma\:\{Q,\rho(g)\}\to\{Q',\rho'(g)\})$ such that
$\tau=\delta(\sigma)$.                                      

To this end, we write $\tau$ as a family of linear maps and translate the fact that $\tau$ is a morphism in $\Cal C$
into conditions in these linear maps.

We write $Q$ and $Q'$ as quivers,
$$ Q=(M_{i},m_{ji}),\quad Q'=( M', m'_{ji}).
$$ The isomorphisms $\rho(g)\:Q\to \Phi(g)Q$ are given by isomorphisms of linear spaces (3) such that the
diagrams (4) commute.

Similarly, for $\{Q',\rho'(g)\}$ we have isomorphisms
$$
\gamma'_{i,g}\:M'_i \to M'_{g(i)}
\tag6
$$ of linear spaces such that for each $i,j,g$ the diagram
$$
\CD M'_i   @>\gamma'_{i,g}>>  M'_{g(i)}\\ @Vm'_{ji}VV   @VV\widetilde m'_{g(j),g(i)}V \\ M'_j
@>\gamma'_{j,g}>> M'_{g(j)}
\endCD
\tag7
$$ obtained from (4) by adding primes, commutes. In particular,  
$\widetilde m'_{g(j),g(i)}$ are expressed in terms of $m'_{kl}$ and $(1+m'_{kk})^{-1}$ by the same polynomials
that give expressions of  $\widetilde m_{g(j),g(i)}$ in terms of $m_{kl}$ and $(1+m_{kk})^{-1}$. 

A morphism $\tau\:\delta(\{Q,\rho(g)\}) \to \delta(\{Q',\rho'(g)\})$ is given by a family of linear maps
$\tau_i\:M_i\to M'_i$ such that the for each $i,g$ the diagram              
           
$$
\CD M_i   @>\gamma_{i,g}>>  M_{g(i)}\\ @V\tau_i VV  @VV\tau_{g(i)}V\\ M'_i   @>\gamma'_{i,g}>> 
M'_{g(i)}
\endCD
\tag8
$$                                    commutes, and for $(i,j)$ the diagram
$$
\CD M_i   @>\tau_i>> M'_i\\ @Vm_{ji} VV  @VVm'_{ji}V\\ M_j   @>\gamma'_{i,g}>>  M'_j
\endCD
\tag9
$$                                    commutes (since $\tau$ is a morphism in $\Cal Q$).

On the other hand, a morphism $\sigma$  is a family of linear maps
$$
\sigma_{i}\:M_{i}\to M'_{i}
$$ such that 

To construct $\sigma\:\{Q,\rho(g)\}\to\{Q',\rho'(g)\}$ we must define  $\sigma_i\: M_i\to M'_i$. We set
$\sigma_i=\tau_i$. We must prove that all cubic diagrams with (4) as the top square, (7) as the bottom square, and
$\sigma$'s as vertical edges (directed down) commute.

In any such cubic diagram, the square corresponding to the top face commutes by (4), the square corresponding to
the bottom face commutes by (7), the squares corresponding to front and back faces commute by (8), and the square
corresponding to the lest face commutes by (9). Since all left-to-right arrows are isomorphisms -- they are various
$\gamma_{i,g}$ --, the square corresponding to the right face also commutes. Therefore, the family
$\sigma_i\:M_i\to M'_i$ indeed determines a morphism
$\sigma\:\{Q,\rho(g)\}\to\{Q',\rho'(g)\}$. Clearly,
$\tau=\delta(\sigma)$. Lemma 5.3 is proved.                  \enddemo

\heading 6. Embedding of microlocal perverse sheaves into meldings \endheading
                                                          
In this section we show that the morphism of stacks $\nu: \ode \to \om$ is an  embedding of stacks. This means in
particular, combined with the results of Section 5, that the stack $\om$ has enough data to describe all microlocal
perverse  sheaves. What remains is to describe the relations coming from the loci of codimension $\geq 2$ in the
support.

We need the following construction due to Kashiwara and Schapira (see \cite{KS}).  Let
$A$ and  $B$ be two sheaves in $\od^b(X)$ and let $i\: \Delta X \hookrightarrow X 
\times X$ be the diagonal embedding. We define $$
\muhom(A,B) = \mu_{\Delta @, X} R\boldkey H \boldkey o \boldkey m  (p_1^{-1}A,p_2^!B),
$$ where $p_i:X \times X \to X$ are the projections. Clearly $\muhom(A,B) \in\od^b(\ct)$ and  
$$
\Hom_{\od ^b(X)}(A,B) \cong H^0(\ct,\muhom(A,B)).
$$ We also have

\proclaim{Proposition 6.1} For $A,B \in \od^b(X)$ we have 
$\text{Supp}(\muhom(A,B)) \subset \oss(A)\cap\oss(B)$
\endproclaim

\demo{Proof} See \cite{KS, 5.4.10 (ii)}.
\enddemo

It is now easy to see that $\muhom$ behaves well with respect to  microlocalization, i.e., for $A,B\in \op(X,U)$,
$U\subset \ct$ open, we have a well defined
$\muhom(A,B) \in \od^b(X,U)$. We also have

\proclaim{Proposition 6.2} If $A,B\in \op(X,U)$, then $\muhom(A,B)[\dim X] \in
\op(U)$.  
\endproclaim

\demo{Proof} The statement being local, we can, by Lemma 2.5, assume that $A$ and
$B$ are perverse sheaves.
 All the constructions used in defining $\muhom$ preserve  the  perversity conditions. For details see \cite{KS,
10.3.20}. 
\enddemo

Given $A,B\in \ode(U)$ we can consider the sheaf $H^0(\muhom(A,B))$ on $U$. By definition of $\muhom$ we
get a morphism of sheaves 
$$
\Hom_{\ode} (A,B) \to H^0(\muhom(A,B)). \tag6.1
$$

\proclaim{Proposition 6.3} The morphism of sheaves in (6.1) is an  isomorphism.
\endproclaim

\demo{Proof} It suffices to prove this statement on stalks. It is then precisely \cite{KS, 6.1.2}.
\enddemo

\proclaim{Theorem 6.4} The morphism $\nu: \ode \to \om$ of stacks is an embedding.
\endproclaim

\demo{Proof} It suffices to prove this statement for a fixed Lagrangian $\Lambda
\subset \ct$. We are then reduced to proving that for $A,B \in  \ode_\Lambda(U)$ the sheaves $\Hom_{\ode}(A,B)$
and 
$\Hom_{\om}(\nu A,\nu B)$ are isomorphic. Let
$Z=\Lambda - (\Lambda^0 \cup \Lambda^1)$, and let $j:U-Z \hookrightarrow U$ be the inclusion. By definition the
sheaf
$\Hom_{\om}(\nu A, \nu B)$ has the property that $\operatorname{Hom}_{\om}(\nu A, \nu  B) \cong  
j_*j^*\Hom_{\om}(\nu A, \nu B)$. By Theorem 5.1, the sheaves
$j^*\Hom_{\ode}(A,B)$ and $j^*\Hom_{\om}(\nu A, \nu B)$ are isomorphic. Therefore, it  suffices to show that
$\Hom_{\ode}(A,B) \cong j_*j^* 
\Hom_{\ode}(A,B)$. By Proposition 6.3, this reduces to the isomorphism 

$$H^0(\muhom(A,B)) \cong j_*j^* 
H^0(\muhom(A,B)). 
$$

To prove this consider the exact sequence
$$ 
\multline H^0_Z(\muhom(A,B))\to H^0(\muhom(A,B)) \to
\\
 j_*j^*H^0(\muhom(A,B)) \to  H^1_Z(\muhom(A,B)).
\endmultline
$$ Since, by proposition 6.2, $\muhom(A,B)[\dim X]$ is perverse, we get that 
$$H^0_Z(\muhom(A,B))=0 \ \text{and} \ H^1_Z(\muhom(A,B))=0\,,$$  
which proves the theorem.
\enddemo

\define\F{\Cal F}

\define\To{\Rightarrow}
\define\Id{\roman {Id}}
\define\id{\roman {id}} 
\define\vv{\varphi}
\vskip 10pt

\heading Appendix on stacks
\endheading
\vskip 10pt

\subhead 0. Introduction \endsubhead In this preface we present a brief account of definitions and results about
stacks of categories that are used in the main text. Roughly speaking, a stack $\Cal F$ on a topological space $X$ is
an analog of a sheaf $F$  on
$X$, with sets $F(U)$ (or sets with additional structure, e.g., groups, or modules) are replaced by categories $\Cal
F(U)$ (or categories with additional structure, e.g., abelian categories). Since the natural notion of ``equality'' of
categories is the notion of equivalence, rather than that of of an isomorphism, one should be careful when
generalizing main constructions and results from sheaves to stacks. 

We do not intend to give here the complete proofs. Instead, the main goal of this appendix is to convince the reader
that, taking some precautions, one can work with stacks as comfortably as with sheaves. For a rigorous and complete
treatment of the subject (although in a somewhat different language) the reader is referred to Chapter 1 in \cite{Gi}.

\definition{1. Definition} A {\it prestack\/} $\F$ on a topological space $X$ is the following collection of data:
\roster
\item"(a)" for each open set $U\subset X$, a category $\F(U)$;

\item"(b)" for each pair of open sets $V\subset U$, a functor
$$
\rho_{VU}\:\F(U)\to \F(V)
$$ (sometimes we think of this functor as restriction to $V$: $A\in \Ob \Cal F(U) \mapsto A\big |_V\in \Ob \F(V)$;

\item"(c)" for any three open sets $W\subset V\subset U$, an isomorphism of functors from $\F(U)$ to $\F(W)$:
$$
\tau_{WVU}\:\rho_{WV}\circ\rho_{VU}\To\rho_{WU}.
$$
\endroster

\noindent These data should satisfy the following conditions: 

\roster
\item"(I.a)" $\rho_{UU}=\Id_{\F(U)}$;

\item"(I.b)" $\tau_{VVU}=\tau_{VUU}=\Id_{\rho_{VU}}$;

\item"(I.c)" for any four open sets $S\subset W\subset V \subset U$ we have $$
\tau_{SWU}\circ\tau_{WVU}=\tau_{SVU}\circ \tau_{SWV}
$$ as morphisms of functors
$$
\rho_{SW}\circ\rho_{WV}\circ\rho_{VU}\To\rho_{SU}.
$$
\endroster\enddefinition

\definition{2. Definition} A {\it morphism of prestacks\/}
$\Theta\:(\F,\rho,\tau)\to (\F',\rho',\tau')$ on a topological space $X$ is a collection of functors $\Theta(U)\:\F(U)\to
\F'(U)$, one for each open set $U\subset X$, and of isomorphisms of functors $\Theta(V,U)\:\Theta(V)\circ
\rho_{VU} \to\rho'_{VU}\circ
\Theta(U)$, such that for any three open sets $W\subset V\subset U$ we have the equality
$$
\tau'_{WVU} \circ \Theta (V,U) \circ \Theta(W,V) = \Theta(W,U)\circ \tau_{WVU} .
$$                                                                            of two isomorphisms of functors
$$
\Theta(W)\circ\rho_{WV}\circ \rho_{VU} \To \rho'_{WU}\circ \Theta(U) $$ (both functors are from the category 
$\F(U)$ to the category
$\F'(W)$). \enddefinition

It is clear how to define the composition of morphisms of prestacks and the identity morphisms, and that with these
definitions prestacks on $X$ form a category.

\definition{3. Definition} A prestack $\Cal F$ on $X$ is called a {\it stack\/} if, in addition to conditions I.a--c of 
Definition 1 it satisfies the following conditions:

\roster

\item"(II)" Let $U\subset X$ be an open set, $A,B$ be two objects of the category
$\F(U)$. Then the correspondence 
$$ U\supset  V\mapsto \Hom_{\F(V)}(\rho_{VU}A, \rho_{VU}B)
$$ is a sheaf on $U$.

\item"(III)" Let $U\subset X$ be an open set, $U=\bigcup U_i$ be an open covering of
$U$. Suppose we are given a family of objects $A_i\in\Ob \F(U_i)$ and of isomorphisms
$$
\sigma_{ij}\:\rho_{U_i\cap U_j,U_i}(A_i) \to \rho_{U_i\cap U_j,U_j}(A_j) $$ (in the category $\F(U_i\cap U_j)$)
such that:

(a) $\sigma_{ii} = \id$;

(b) $\sigma_{ji}= (\sigma_{ij})^{-1}$;

(c) for any three open sets $U_i,U_j,U_k$ and three objects $A_i\in\Ob \F(U_i)$, 
$A_j\in\Ob \F(U_j)$,  $A_k\in\Ob \F(U_k)$, the restrictions of isomorphisms $\sigma_{ij}$, $\sigma_{jk}$, 
$\sigma_{ki}$ to $U_i\cap U_j \cap U_k$ satisfy the condition
$$
\rho_{U_i\cap U_j \cap U_k, U_i\cap U_j}(\sigma_{ij}) \circ
\rho_{U_i\cap U_j \cap U_k, U_j\cap U_k}(\sigma_{jk}) \circ
\rho_{U_i\cap U_j \cap U_k, U_k\cap U_i}(\sigma_{ki}) = \id
$$ (the identity morphism of the object $\rho_{U_i\cap U_j \cap U_k, U_i\cap U_j} \circ
\rho_{U_i\cap U_j, U_i}(A_i)$ in the category $\F(U_i\cap U_j \cap U_k)$).
 
 Then there exists an object $A\in\Ob \F(U)$ and a family of isomorphisms  $$
\sigma_i\: \rho_{U_i,U}(A) \to A_i
$$ in $\F(U_i)$ that are compatible with $\sigma_{ij}$ in the following sense: For any $i,j$ the following diagram
of isomorphisms in the category $\F(U_i\cap U_j)$ is commutative:
$$
\CD
\rho_{U_i\cap U_j,U_i}\circ \rho_{U_i,U} (A)  @> \rho_{U_i\cap U_j,U_i}(\sigma_i) >{\cong}>  \rho_{U_i\cap
U_j,U_i}(A_i) \\ @VV
\tau_{U_i\cap U_j,U_i,U}(A)V  @VVV\\
\rho_{U_i\cap U_j,U}(A)  @. \qquad{ \ \ \ \sigma_{ij}}\\ @AA  {\tau_{U_i\cap U_j,U_j,U}(A)}
A  @VVV\\
\rho_{U_i\cap U_j,U_j}\circ \rho_{U_j,U} (A)  @> \rho_{U_i\cap U_j,U_j}(\sigma_j) >{\cong}>  \rho_{U_i\cap
U_j,U_j}(A_j)
\endCD
$$

The object $A$ is unique up to a unique isomorphism. 
\endroster

Defining morphisms of stacks as morphisms of corresponding prestacks, we obtain the category $\Cal S(X)$ of
stacks on a topological space $X$.

\enddefinition

\definition{4. Definition}  A morphism of prestacks $\Theta\:(\F,\rho,\tau)\to (\F',\rho',\tau')$ is called a {\it weak
isomorphism\/} if all the functors $\Theta(U)$ are equivalences of the corresponding categories.
\enddefinition

Denote the family of all weak isomorphisms of prestacks  by $WI$.

The next result shows that $WI$ is a reasonable family to form the corresponding calculus of fractions.

\proclaim{Proposition} $WI$ is a localizing system of morphisms in the category $\Cal S(X)$ {\rm (see
\cite{GM})}.
\endproclaim

\demo{Sketch of the proof} The main condition to be verified is the Ore condition that allows us to replace left
fractions with right fractions. Namely, we must prove that for two composable morphisms $\Theta,\Psi$ in $\Cal
S(X)$ with $\Psi\in WI$ there exist composable morphisms $\Theta',\Psi'$ in $\Cal S(X)$ with $\Psi'\in WI$ such
that
$\Theta\Psi = \Psi'\Theta'$.

The main part of the proof of the Ore condition is the proof of the corresponding statement for categories instead of
stacks. Namely, given three categories $\Cal A, \Cal B, \Cal C$ and two functors $F\:\Cal A\to \Cal B$, $G\:\Cal
C\to \Cal B$ such that $F$ is an equivalence of categories, we want to construct a category $\Cal B'$ and functors
$G'\:\Cal B'\to \Cal A$, $F'\:\Cal B'\to \Cal C$ such that $F'$ is an equivalence of categories and the diagram $$
\CD
\Cal B'  @>F'>>  \Cal C \\ @VG'VV     @VGVV \\
\Cal A  @>F>>  \Cal B
\endCD
\tag A1
$$ commutes.

We construct $\Cal B'$, $F'$, $G'$ as follows. An object of $\Cal B'$ is a triple
$(A,C,\vv)$, where $A\in \Ob \Cal A$, $C\in \Ob\Cal C$, $\vv\:F(A)\to G(C)$ an isomorphism in $\Cal B$. A
morphism $(A,C,\vv)\to (A_1,C_1,\vv_1)$ in $\Cal B'$ is a pair of morphisms $f\:A\to A_1$ in $\Cal A$, $g\:C\to
C_1$ in $\Cal C$ such that the diagram
$$
\CD F(A)  @>F(f)>>  F(A_1)\\ @V\vv VV   @V\vv_1 VV\\ G(C)  @>G(g)>>  G(C_1)
\endCD
$$ commutes.

Functors $G'\:\Cal B'\to \Cal A$ and $F'\:\Cal B'\to \Cal C$ are defined by 
$(A,C,\vv)\mapsto A$ and $(A,C,\vv)\mapsto C$ respectively. The commutativity of (A1) is clear.To prove that 
$F'$ is an equivalence of categories, we  apply the Freyd theorem (see \cite{M}), which says that a functor $K\:\Cal
D_1\to \Cal D_2$ between two categories is an equivalence of  categories if and only if it induces isomorphisms of
Hom- sets and each object from $\Cal D_2$ is isomorphic to an object of the form $K(D_1)$ for some object $D_1$
from $\Cal D_1$. The corresponding properties of the functor $F'$ from diagram (A1) immediately follow from the
fact that $F$ is an equivalence of categories.

Since we have constructed $\Cal B'$, $G'$, and $F'$ canonically, it is easy to generalize this construction to stacks,
thus proving the Ore condition for the class $WI$.
\enddemo

Denote by $\widetilde{\Cal S}(X)$ the localized category  
$\widetilde{\Cal S}(X) = \Cal S(X)[WI]^{-1}$.

\subhead 5. Stackification of a prestack \endsubhead Similarly to sheaves, given a prestack $\Cal F$ on a topological
space $X$, one can construct a canonical stack $\Cal G$ associated to $\Cal F$. 

\proclaim{Theorem} Let $\Cal F$ be a prestack on $X$. Then there exists a stack  $\Cal G$ on $X$ and a morphism
of prestacks
$\alpha\:\Cal F\to \Cal G$ which is universal in the following sense: For any stack $\Cal E$ on $X$ and any
morphism of prestacks
$\varphi\:\Cal F\to \Cal E$ there exists a unique morphism of stacks $\psi\:\Cal G\to
\Cal E$ such that $\varphi=\psi \circ\alpha$. If $(\Cal G,\alpha)$, $(\Cal G',\alpha')$ are two such pairs, then there
exist unique  weak isomorphisms $\Cal G\to \Cal G'$,
$\Cal G'\to \Cal G$ commuting with $\alpha$, $\alpha'$.  
\endproclaim

We call $\Cal G$ the {\it stack associated to the prestack $\Cal F$}. If $\Cal F$ is itself a stack, then clearly $\Cal
G=\Cal F$ and
$\alpha=\id$.

The above theorem is one of the reasons why the localized category $\widetilde{\Cal S}(X)$ is more convenient
than the category $\Cal S(X)$: in $\widetilde{\Cal S}(X)$ the stackifization is unique up to a unique isomorphism.

\subhead 6. The Locality theorem \endsubhead The next theorem shows that stacks are analogous to sheaves in the
sense that when we need to perform certain  constructions and/or verify certain properties  of stacks, it suffices to do
this on small open sets.

We say that we are given a local morphism of stacks $\F\to \F'$ if the functors
$\Theta(U)$ in Definition 2 are given only for small open sets $U$. More precisely, we assume that for each point
$x\in X$ we have a neighborhood $U_x$ and that
$\Theta(U)$, $\Theta(V,U)$ are defined for open sets $U,V$ such that $V\subset U\subset U_x$ for some $x$, and
satisfy the conditions of Definition 2 for such open sets.

\proclaim{Locality Theorem} Any local morphism $\Theta\:\F\to \F'$ of stacks can be uniquely extended to a
morphism
$\overline\Theta\:\F\to \F'$ of stacks. If $\Theta$ is a weak isomorphism, then $\overline\Theta$ is also a weak
isomorphism.  
\endproclaim

\demo{Sketch of the proof} 1. To construct the extension of a local morphism of stacks to a morphism of stacks we
must define
$\Theta(U)$ and $\Theta(V,U)$ for larger open sets. Assume that $\Theta(U_1)\:\F(U_1)\to \F'(U_1)$ and
$\Theta(U_2)\:\F(U_2)\to
\F'(U_2)$ are already defined. Denote $U=U_1\cup U_2$, $\widehat {U}=U_1\cap U_2$, and define the functor
$\Theta(U)\:\F(U)\to
\F'(U)$.

Let $A\in \Ob \F(U)$. Denote $A_1=\rho_{U_1,U}(A) \in \Ob\F(U_1)$,    
$A_2=\rho_{U_2,U}(A) \in \Ob\F(U_2)$, and let the isomorphism
$$
\sigma\: \rho_{\widehat {U},U_1}(A_1) \to  \rho_{\widehat {U},U_2}(A_2) $$ be given by the formula
$$
\sigma = (\tau_{\widehat {U} U_2U})^{-1}\circ \tau_{\widehat {U} U_1U} $$ Let $A'_1 = \Theta(U_1)(A_1)$,
$A'_2 = \Theta(U_2)(A_2)$. There is a natural way to construct an isomorphism
$$
\sigma'\: \rho'_{\widehat {U},U_1}(A'_1) \to  \rho'_{\widehat {U},U_2}(A'_2) $$ as the compositions of
isomorphisms $\sigma$,
$\Theta(\widehat {U},U_1)$, 
$\Theta(\widehat {U},U_2)$, and their inverses (it is easier for the reader to figure out how to do it than to read
complicated formulas).

By condition (III) of Definition 3 applied to the stack $\F'$, the triple $(A'_1,A'_2,\sigma')$ yields a unique object
$A'\in
\Ob\F'(U)$, and we take $A' = \Theta(U)(A)$. The possibility to extend $\Theta(U)$ to a functor $\F(U)\to \F'(U)$
and the construction of functor isomorphisms $\Theta(V,U)$ for $V\subset U$ is, of course, left to the reader.

2. To prove the second statement of the theorem we must prove that if both
$\Theta(U_1)$ and $\Theta(U_2)$ are equivalences of categories, then the functor
$\Theta(U)$ is also an equivalence of categories. By the Freyd theorem, it suffices to prove that $\Theta(U)$ is a
bijection on
$\Hom$-sets and that any object $A'$ in $\F'(U)$ is isomorphic to an object  of the  form $\Theta(U)(A)$ for some
$A\in\Ob \F(U)$. The first statement follows from Condition (II) of Definition 3. The second statement can be
proved similarly to the construction of the functor
$\Theta(U)$ above.  
\enddemo

\subhead 7. The stalk of a stack \endsubhead Let $\F\in \Cal S(X)$,  $x$ a point in $X$. The {\it stalk\/} of $\F$ at
$x$ is the category $\F_x$ that is defined as follows. An object of $\F_x$ is a pair $(U,A)$, where $U\subset X$ is
an open set containing $x$,
$A\in \Ob\F(U)$. A morphism $(U,A)\to (V,B)$ is an equivalence class of pairs
$(R,\alpha)$, where $x\in R\subset U\cap V$, $\alpha\:A\big |_R\to B\big |_R$ a morphism in $\F(R)$, with two
pairs $(R,\alpha)$ and
$(S,\beta)$ being equivalent if there exists a third pair $(T,\gamma)$ such that $x\in T\subset R\cap S$ and
$\gamma=\alpha\big | T = \beta\big |_T$. The composition
$$ (U,A) @>(R,\alpha)>>  (V,B)  @>(S, \beta)>> (W,C)
$$ is the class of the pair $(T,\gamma)$, where $T$ is an arbitrary open set with $x\in T\subset R\cap S$ and
$\gamma = \beta\big | _T\circ \alpha\big | _T$. Clearly, different choices of $T$ give equivalent pairs $(T,\gamma)$.

The identity morphism $(U,A)\to (U,A)$ is the class of the pair $(U,\id_A)$.

A morphism of stacks $\Theta\:\F\to \F'$ determines in an evident way a functor
$\Theta_x\:\F_x\to \F'_x$.

\proclaim{8. Theorem} Let $\Theta\:\F\to \F'$ be a weak isomorphism of stacks. Then, for any $x\in X$, $\Theta_x$
is an equivalence of categories. \endproclaim

\demo{Proof} (a) Let $(U,A')\in\Ob \F'_x$. Since $\Theta(U)\:\F(U)\to \F'(U)$ is an equivalence of categories, $A'$
is isomorphic in
$\F'(U)$ to an object of the form $\Theta(U)(A)$ for some $A\in \Ob\F(U)$. Then
$(U,A')$ is isomorphic in $\F'_x$ to $\Theta_x(U,A)$.

(b) Let $(U,A),(V,B)\in \Ob\F_x$, $A'=\Theta(U)A$, $B'=\Theta(V)B$, so that
$(U,A'),\allowmathbreak(V,B') \in \Ob\F'_x$ and $(U,A')=\Theta_x(U,A)$, $(V,B') =\Theta_x(V,B)$.

Let $\vv'\:(U,A')\to (V,B')$ be a morphism in $\F'_x$. It is represented by a morphism
$\overline{\vv}{}'\:A'\big |_S\to B'\big |_S$ in $\F'(S)$ for some   $S\subset U\cap V$, $x\in S$. Since $\Theta(S)$
is an equivalence of categories, there exists a unique morphism $\overline{\vv}\:A\big |_S\to B\big |_S$ in $\F(S)$
such that
$\overline{\vv}{}' = \Theta(S)\overline{\vv}$. This $\overline{\vv}$ yields a morphism
$\vv\:(U,A)\to (V,B)$ in $\F_x$ such that $\vv' = \Theta_x \vv$. 

(c) Let us show that the morphism constructed in (b) is unique. Indeed, another such morphism $\psi$ is given by
some $T\subset U\cap V$, $x\in T$, and a morphism
$\overline{\psi}\:A\big |_T\to B\big |_T$ such that the class of pair
$(T,\Theta(T)\overline{\psi})$ is the same morphism $\vv'\:(U,A')\to(V,B')$ in $\F'_x$. This means that for some
open set $W\subset S\cap T$ we have $\overline{\vv}{}'\big |_W = \Theta(T)\overline{\psi}\big |_W$ in the
category $\F'(W)$. Since $\Theta(W)$ is an equivalence, 
$$
\overline{\vv}\big |_W = \overline{\psi}\big |_W,
$$ i.e., the class of pair $(T,\psi)$ in $\F_x$ coincides with $\vv$.

By the Freyd theorem, (a)--(c) imply that $\Theta_x\:\F_x\to \F'_x$ is an equivalence of categories.
\enddemo

Now we prove the converse result.

\proclaim{9. Theorem} Let $\Theta\:\F\to \F'$ be a morphism of stacks on a topological space $X$ such that for any
$x\in X$ the functor $\Theta_x:\F_x\to \F'_x$ is an equivalence of categories. Then $\Theta$ is a weak equivalence.
\endproclaim

\demo{Proof} It suffices to prove that $\Theta(X)\:\F(X)\to \F'(X)$ is an equivalence of categories. Again, we use
the Freyd theorem. Let $A,B\in \Ob\Cal \F(X)$,
$A'=\Theta(X)(A)$, $B'=\Theta(X)(B)$, and $\vv'\:A'\to B'$ a morphism in $\F(X)$.  For any $x\in X$, denote
$(X,A)=A_x$, $(X,B)=B_x$ and similarly for $(X,A')$, $(X,B')$. The morphism $\vv'(x)\: A'_x\to B'_x$
corresponding to $\vv'$ comes, via $\Theta_x$, from a unique morphism  $\vv(x)\: A_x\to B_x$. This means that
there exists a covering
$\{U_i\}$ of $X$ by open sets and morphisms $\vv_i\:A\big |_{U_i} \to B\big |_{U_i}$ such that
$\Theta(\vv_i)=\vv'\big |_{U_i}$. 

Similarly to the proof of the locality theorem above, we use Condition (II) in Definition 3 to conclude that the
morphisms $\vv_i$ glue together to a unique morphism $\vv\:A\to B$ such that $\Theta(X)(\vv)=\vv'$. Therefore is
remains to prove that each $A'\in \Ob
\F'(X)$ is isomorphic to an object of the form $\Theta(A)$ for $A\in \Ob \F(X)$. Choose an arbitrary $x\in X$.
Since $\Theta_x$ is an equivalence of categories for each $x\in X$, there exist an open covering $\{U_i\}$ of $X$,
objects $A_i\in\Ob \F(U_i)$, and isomorphisms $\vv_i\:\Theta(A_i)\to A'\big |_{U_i}$ in $\F'(U_i)$. For $i,j$ with
nonempty intersection $U_{ij}=U_i\cap U_j$ denote by $\sigma_{ij}$ the composite isomorphism
$$ A_i\big | _{U_{ij}} @>\vv_i\big |_{U_{ij}}>>  A' \big | _{U_{i}} \big | _{U_{ij}}  @> \tau_{XU_iU_{ij}}>> 
A'\big | _{U_{ij}}  @>(\tau_{XU_jU_{ij}})^{-1}>>  A'\big | _{U_{j}}\big | _{U_{ij}}  @>(\vv_j\big
|_{U_{ij}})^{-1}>> A_j\big | _{U_{ij}}
$$ The family $(A_i,\sigma_{ij})$ satisfies Condition (III) of Definition 3. Since $\F$ is a stack, there exist
$A\in\Ob \F$ and isomorphisms $\sigma_i\:A \big |_{U_{i}} \to A_i$ that are compatible with $\sigma_{ij}$ in the
sense of this Condition (III). One can easily see that $\Theta(X)(A)=A' $.
\enddemo

\subhead 10. Direct and inverse images \endsubhead Similarly to sheaves, one can define direct and inverse images
of stacks under continuous map of topological spaces. General definition requires a lot of technical details.
Fortunately, in the present paper we need only inverse image  of a stack under open maps, and this is easy.

\definition{Definition} Let $f\:X\to Y$ be a continuous map of topological spaces such that for any open set
$U\subset Y$ its image $f(U)$ is open in $Y$. Let $\Cal F$ be a stack on $Y$. Then its inverse image $f^*\Cal F$ is
defined by the formula $(f^*\Cal F)(U) = \Cal F(f(U))$.

It is clear that the categories $f^*\Cal F(U)$, together with evident restriction functors
$\rho_{VU}$, form a stack on $X$.
\enddefinition

An example of a map satisfying the conditions of the previous definition (the only case we need in this paper) is the
projection of the total space of a locally trivial bundle to the base.

\enddocument